\documentclass[12pt]{article}
\usepackage{amsmath,amssymb,amsthm}
\usepackage{hyperref}
\usepackage{datetime}
\usepackage{units}
\usepackage{color}
\usepackage[T1]{fontenc}
\usepackage[utf8]{inputenc}
\usepackage{authblk}
\usepackage{bm,latexsym,mathrsfs,enumerate}
\usepackage{color}

\title{ A novel approach to the discovery of binary BBP-type formulas for polylogarithm constants\thanks{%
MSC 2010: 11Y60, 30B99}}
\author[]{Kunle Adegoke\thanks{adegoke00@gmail.com }}
\affil{Department of Physics and Engineering Physics, \mbox{Obafemi Awolowo University, Ile-Ife, 220005 Nigeria}}

\theoremstyle{plain}
\numberwithin{equation}{section}

\begin{document}
\date{}
\maketitle
\begin{abstract}
\noindent Using a clear and straightforward approach, we discover and prove new binary digit extraction BBP-type formulas for polylogarithm constants. Some known results are also rediscovered in a more direct and elegant manner. Numerous experimentally discovered and previously unproved binary BBP-type formulas are also proved. 
\end{abstract}
\tableofcontents

\section{Introduction}	

This paper is concerned with proving each of a lengthy list of conjectured binary BBP-type formulas collected in the ``Compendium of BBP-Type Formulas'', an online collection of BBP-type formulas for various mathematical constants~\cite{bailey09}. The formulas have, in most cases, been outstanding for upwards of fifteen years, in spite of many thousands of downloads of the BBP Compendium. New binary BBP-type formulas, together with their proofs will also be derived.

\bigskip

BBP-type formulas are formulas of the form
\[
c = \sum_{k=0}^\infty  1/b^k \sum_{j=1}^l a_j / (k l + j)^s
\]
where $s$, $b$, $l$ and $a_j$ are integers, and $c$ is some constant. Formulas of this type were first introduced in a 1996 paper~\cite{bbp97}, where a formula of this type for $\pi$ was given. Such formulas allow digit extraction --- the $i$-th digit of a mathematical constant $c$ in base $b$ can be calculated directly, without needing to compute any of the previous $i-1$ digits, by means of simple algorithms that do not require multiple-precision arithmetic~\cite{bailey09}.

\bigskip

Apart from digit extraction, another reason the study of BBP-Type formulas has continued to attract attention is that BBP-Type constants are conjectured to be either rational or normal to base $b$~\cite{bailey01,borwein02,chamberland03}, that is their base-$b$ digits are randomly distributed. 

\bigskip

BBP-Type formulas are usually discovered experimentally, through computer searches, by using Bailey and Ferguson's PSLQ (Partial Sum of Squares -- Lower Quadrature) algorithm~\cite{ferguson99} or its variations. PSLQ and other integer relation finding schemes typically do not suggest proofs~\cite{borwein02,bailey06}. Formal proofs must be developed after the formulas have been discovered.

\bigskip

Our approach in this paper is the systematic symbolic (that is, non-computer-search-based) discovery of BBP-type formulas. The methods used here aim to complement the experimental approaches that have dominated the area. Through fundamental methods, a wide range of interesting, mostly new, BBP-type formulas will be obtained, together with their proofs. It should be noted that in this paper, unlike in \cite{bbp97} or \cite{huvent01b}, for example, no evaluation of complicated integrals is necessary. The BBP-type formulas come as natural consequences of the corresponding polylogarithm identities.

\bigskip
\section{Notation}
Degree~$s$ ($\in\mathbb{Z^+}$in this paper) polylogarithm functions ${\rm Li}_s$ are defined by
\[
{\rm Li}_s [z] = \sum\limits_{k = 1}^\infty  {\frac{{z^k }}{{k^s }}},\quad |z|\le 1\,. 
\]

In particular, for $|z|=1$ and $x\in\mathbb{R}$ we have
\begin{equation}
\begin{split}
{\rm Li}_{2n} [e^{ix} ]& = {\rm Gl}_{2n} (x) + i{\rm Cl}_{2n} (x)\\
{\rm Li}_{2n + 1} [e^{ix} ]& = {\rm Cl}_{2n + 1} (x) + i{\rm Gl}_{2n + 1} (x)\,,
\end{split}
\end{equation}
where ${\rm Gl}$ and ${\rm Cl}$ are Clausen sums~\cite{lewin81} defined, for $n\in\mathbb{Z^+}$ by

\begin{equation}
\begin{split}
{\rm Cl}_{2n} (x)& = \sum\limits_{k = 1}^\infty  {\frac{{\sin kx}}{{k^{2n} }}},\quad {\rm Cl}_{2n + 1} (x) = \sum\limits_{k = 1}^\infty  {\frac{{\cos kx}}{{k^{2n + 1} }}}\\
{\rm Gl}_{2n} (x) &= \sum\limits_{k = 1}^\infty  {\frac{{\cos kx}}{{k^{2n} }}},\quad {\rm Gl}_{2n + 1} (x) = \sum\limits_{k = 1}^\infty  {\frac{{\sin kx}}{{k^{2n + 1} }}}\,.
\end{split}
\end{equation}

We shall find the following formulas useful:

\begin{equation}
\begin{split}
{\rm Gl}_{2n} (x) &= ( - 1)^{1 + [n/2]} 2^{n - 1} \pi ^n {\rm B}_n (x/2\pi )/n!\\
\frac{1}{{m^{n - 1} }}{\rm Cl}_n (mx) &= \sum\limits_{r = 0}^{m - 1} {{\rm Cl}_{n} (x + 2\pi r/m)}\,.
\end{split}
\end{equation}
Here $[n/2]$ denotes the integer part of $n/2$ and $ {\rm B}_n$ are the Bernoulli polynomials defined by

\[
{\frac {t{{\rm e}^{xt}}}{{{\rm e}^{t}}-1}}=\sum _{n=0}^{\infty }{
\frac {{\rm B}_n (x) {t}^{n}}{n!}}\,.
\]

\bigskip

In order to save space, we will give the BBP-type formulas using the compact P-notation~\cite{bailey09}:

\[
\sum\limits_{k = 0}^\infty  {\frac{1}{{b^k }}\sum\limits_{j = 1}^l {\frac{{a_j }}{{(kl + j)^s }}} }\equiv P(s,b,l,A) \,, 
\]
where $s$, $b$ and $l$ are integers, and \mbox{$A = (a_1, a_2,\ldots, a_l)$} is a vector of integers.

\section{Scheme for obtaining the BBP-type formulas}\label{sec.scheme}

The derivation of a desired degree $s$ BBP-type formula proceeds in two stages

\begin{enumerate}
\item An attempt is made to express the polylogarithm constant, $c$, of interest as a linear combination of the real or imaginary parts of polylogarithms:

\begin{equation}
c = \sum\limits_j {\left\{ {\alpha _j\, {\rm Re\;Li}_s \left[ {p_j \exp (ix_j )} \right]} \right\}} 
\end{equation}

or

\begin{equation}
c = \sum\limits_j {\left\{ {\beta _j\, {\rm Im\;Li}_s \left[ {q_j \exp (iy_j )} \right]} \right\}} 
\end{equation}

for $\alpha _j$, $\beta _j\in\mathbb{Q}$, $p_j$, $q _j\in (0,1)$ and $x_j$, $y_j$ rational multiples of $\pi$.

\item The identities 

\begin{equation}\label{equ.jvdhotv}
{\mathop{\rm Re\;}\nolimits} {\rm Li}_s \left[ {pe^{ix} } \right] = \sum\limits_{k = 1}^\infty  {\frac{{p^k \cos kx}}{{k^s }}}
\end{equation}
and
\begin{equation}\label{equ.rvyvo35}
{\mathop{\rm Im\;}\nolimits} {\rm Li}_s \left[ {pe^{ix} } \right] = \sum\limits_{k = 1}^\infty  {\frac{{p^k \sin kx}}{{k^s }}}\,,
\end{equation}
for $p\in [0,1]$, $x\in [0,2\pi]$ and $s\in\mathbb{Z^+}$, for the real and imaginary parts of a polylogarithm function are then employed to write each constituent term of the linear combination as a BBP-type formula and the indicated combination is then formed. In particular if a binary formula is sought, then $p_j$ and $q_j$ in the above formulas must be taken as positive integral powers of $1/2$ for $x=\pi$, $x=\pi/2$ or $x=\pi/3$ (odd positive integral powers of $1/\sqrt 2$ for $x=\pi/4$ or $x=3\pi/4$).

\end{enumerate}

\bigskip

To accomplish the first stage, polylogarithm functional equations are evaluated at certain carefully chosen coordinates and real and imaginary parts are taken. Sometimes it will be necessary to simultaneously solve two or more equations that couple several polylogarithm constants.

\bigskip

As a concrete example of how Eqs.~\eqref{equ.jvdhotv} and \eqref{equ.rvyvo35} give rise to BBP-type formulas, consider the choice $p=1/\sqrt 2^q$ and $x=\pi/4$ in Eq.~\eqref{equ.jvdhotv} for $q\in\mathbb{Z^+}$ and $\mod(q,2)=1$. This choice, together with the periodicity of $\cos(\pi/4)$ yield the very general BBP-type formula (with $q$ a positive odd integer):

\begin{eqnarray}\label{equ.kyh42ro}
 {\rm Re\,Li}_s \left[ {\frac{1}{{\sqrt 2 ^q }}\exp \left( {\frac{{i\pi }}{4}} \right)} \right]& =& \frac{1}{{2^{12q - s} }}P(s,2^{12q} ,24,(2^{ - {\textstyle{1 \over 2}} - {\textstyle{q \over 2}} + 12q} ,0, - 2^{ - {\textstyle{1 \over 2}} - {\textstyle{q \over 2}} + 11q},\nonumber\\ 
&&  - 2^{10q} , - 2^{ - {\textstyle{1 \over 2}} - {\textstyle{q \over 2}} + 10q} ,0,2^{ - {\textstyle{1 \over 2}} - {\textstyle{q \over 2}} + 9q} ,2^{8q} ,2^{ - {\textstyle{1 \over 2}} - {\textstyle{q \over 2}} + 8q} ,0, \nonumber\\ 
&&  - 2^{ - {\textstyle{1 \over 2}} - {\textstyle{q \over 2}} + 7q} , - 2^{6q} , - 2^{ - {\textstyle{1 \over 2}} - {\textstyle{q \over 2}} + 6q} ,0,2^{ - {\textstyle{1 \over 2}} - {\textstyle{q \over 2}} + 5q} ,2^{4q} , \nonumber\\ 
&&  2^{ - {\textstyle{1 \over 2}} - {\textstyle{q \over 2}} + 4q} ,0,- 2^{ - {\textstyle{1 \over 2}} - {\textstyle{q \over 2}} + 3q} , - 2^{2q} , - 2^{ - {\textstyle{1 \over 2}} - {\textstyle{q \over 2}} + 2q} ,0,\nonumber\\
&& 2^{ - {\textstyle{1 \over 2}} + {\textstyle{q \over 2}}} ,1))\,. 
 \end{eqnarray}

We also have (with $x=3\pi/4$)

\begin{eqnarray}\label{equ.fczffj2}
 {\rm Re\;Li}_s \left[ {\frac{1}{{\sqrt 2 ^q }}\exp \left( {\frac{{3i\pi }}{4}} \right)} \right] &=& \frac{1}{{2^{12q - s} }}P(s,2^{12q} ,24,( - 2^{ - {\textstyle{1 \over 2}} - {\textstyle{q \over 2}} + 12q} ,0,2^{ - {\textstyle{1 \over 2}} - {\textstyle{q \over 2}} + 11q} ,\nonumber\\ 
&& - 2^{10q} ,2^{ - {\textstyle{1 \over 2}} - {\textstyle{q \over 2}} + 10q} ,0, - 2^{ - {\textstyle{1 \over 2}} - {\textstyle{q \over 2}} + 9q} ,2^{8q} , - 2^{ - {\textstyle{1 \over 2}} - {\textstyle{q \over 2}} + 8q} ,0,\nonumber\\ 
&& 2^{ - {\textstyle{1 \over 2}} - {\textstyle{q \over 2}} + 7q} , - 2^{6q} ,2^{ - {\textstyle{1 \over 2}} - {\textstyle{q \over 2}} + 6q} ,0, - 2^{ - {\textstyle{1 \over 2}} - {\textstyle{q \over 2}} + 5q} ,2^{4q} ,\nonumber\\ 
&&- 2^{ - {\textstyle{1 \over 2}} - {\textstyle{q \over 2}} + 4q} ,0,2^{ - {\textstyle{1 \over 2}} - {\textstyle{q \over 2}} + 3q} , - 2^{2q} ,2^{ - {\textstyle{1 \over 2}} - {\textstyle{q \over 2}} + 2q} ,0,\nonumber\\ 
&&- 2^{ - {\textstyle{1 \over 2}} + {\textstyle{q \over 2}}} ,1))\,. 
 \end{eqnarray}

Another example is ($q\in\mathbb{Z^+}$)

\begin{align}\label{equ.r0xlgws}
{\rm Li}_s \left[ { - \frac{1}{{2^q }}} \right] &= {\rm Re\;Li}_s \left[ {\frac{1}{{2^q }}\exp \left( {i\pi } \right)} \right] \nonumber\\ 
& = \frac{1}{{2^{12q - s} }}P(s,2^{12q} ,24,(0, - 2^{11q} ,0,2^{10q} ,0, - 2^{9q} ,0,2^{8q} ,0, - 2^{7q} ,\nonumber\\ 
&\qquad 0,2^{6q} ,0, - 2^{5q} ,0,2^{4q} ,0, - 2^{3q} ,0,2^{2q} ,0, - 2^q ,0,1))\,. 
\end{align}

The reader should note that the series given above are not the only possible ones for the indicated constants; in general, the series used will depend on the particular base and length that are targeted. For instance, a base $2^{4q}$, length $8$ version of \eqref{equ.fczffj2}, for $q$ an odd positive integer, is

\begin{eqnarray}\label{equ.ibbh6bz}
 {\rm Re\;Li}_s \left[ {\frac{1}{{\sqrt 2 ^q }}\exp \left( {\frac{{3i\pi }}{4}} \right)} \right] &=& \frac{1}{{2^{4q} }}P(s,2^{4q} ,8,( - 2^{ - {\textstyle{1 \over 2}} - {\textstyle{q \over 2}} + 4q} ,0,2^{ - {\textstyle{1 \over 2}} - {\textstyle{q \over 2}} + 3q} ,\nonumber \\ 
&& - 2^{2q} ,2^{ - {\textstyle{1 \over 2}} - {\textstyle{q \over 2}} + 2q} ,0, - 2^{ - {\textstyle{1 \over 2}} + {\textstyle{q \over 2}}} ,1))\,,
\end{eqnarray}

while a base $2^{4q}$, length $24$ version of the same series is

\begin{eqnarray}
 {\rm Re\;Li}_s \left[ {\frac{1}{{\sqrt 2 ^q }}\exp \left( {\frac{{3i\pi }}{4}} \right)} \right] &=& \frac{{3^s }}{{2^{4q} }}P(s,2^{4q} ,24,(0,0, - 2^{ - {\textstyle{1 \over 2}} - {\textstyle{q \over 2}} + 4q} ,0,0,0,0,0,\nonumber \\ 
&& 2^{ - {\textstyle{1 \over 2}} - {\textstyle{q \over 2}} + 3q} ,0,0, - 2^{2q} ,0,0,2^{ - {\textstyle{1 \over 2}} - {\textstyle{q \over 2}} + 2q} ,0,0,0,0,0,\nonumber\\
&& - 2^{ - {\textstyle{1 \over 2}} + {\textstyle{q \over 2}}} ,0,0,1))\,. 
\end{eqnarray}

It is of course possible to give BBP-type formulas in general bases for other classes of polylogarithm constants. This is however not the subject matter of this paper.

\bigskip

By way of a specific illustration of how to derive a BBP-type formula for a polylogarithm constant, let us apply the above procedure to obtain a base $2^{12}$ length $24$ formula for $\log^2 2$. The first step is to express this constant as a linear combination of polylogarithms. This is accomplished through the identity (Eq.~\eqref{equ.c5ut706}, subsection~\ref{sec.gendeg2})

\[
\begin{split}
\log ^2 2 &= 2\,{\rm Li}_2 \left[ { - \frac{1}{4}} \right] - 4\,{\rm Re\; Li}_{\rm 2} \left[ {\frac{1}{{\sqrt 2 }}\exp \left( {\frac{{3\pi i}}{4}} \right)} \right]\\
&\qquad- 4\,{\rm Re\; Li}_{\rm 2} \left[ {\frac{1}{{\left( {\sqrt 2 } \right)^3 }}\exp \left( {\frac{{\pi i}}{4}} \right)} \right]\,.
\end{split}
\]

The next step is to now write each of the three constituent members on the right hand side as a BBP-type formula and then form the indicated combination. The result is

\[
\begin{split}
\log^2 2 &=\frac{1}{2^{10}}P(2,2^{12},24,(2^{11},0,-5\cdot2^{11},-7\cdot2^{10},-2^{9},0,2^{8},7\cdot2^{8},5\cdot2^{8},0,-2^{6},2^{7},-2^{5},0,\\
&\qquad 5\cdot2^{5},7\cdot2^{4},2^{3},0,-2^{2},-7\cdot2^{2},-5\cdot2^{2},0,1,-2))\,.
\end{split}
\]

\section{Degree $1$ Formulas}

Degree $1$ BBP-type formulas in general bases are discussed in~\cite{adegokebbp2}. Binary formulas are easily obtained by choosing bases that are powers of $2$. 
Degree $1$ formulas will not be discussed further in this paper.

\section{Degree $2$ Formulas}

\subsection{Generators of Degree~$2$ BBP-type Formulas}\label{sec.gendeg2}

The dilogarithm reflection formula (\mbox{Eq.~A.2.1.7} of~\cite{lewin81}) is 

\[
\frac{{\pi ^2 }}{6} - \log x\log (1 - x)={\rm Li}_2 [x] + {\rm Li}_2 [1 - x]\,.
\]

Putting $x=1/2$ in the above formula gives the well-known result:
\begin{equation}\label{equ.l71m5vv}
\frac{{\pi ^2 }}{{12}} - \frac{{\log ^2 2}}{2} = {\rm Li}_2 \left[ {\frac{1}{2}} \right]
\end{equation}

A two-variable functional equation for dilogarithms, due to Kummer (\mbox{Eq.~A.2.1.19} of~\cite{lewin81}) is

\begin{equation}\label{equ.zbzsvtl}
\begin{split}
{\rm Li}_2 \left[ {\frac{{x(1 - y)^2 }}{{y(1 - x)^2 }}} \right] &= {\rm Li}_2 \left[ { - \frac{{x(1 - y)}}{{(1 - x)}}} \right] + {\rm Li}_2 \left[ { - \frac{{(1 - y)}}{{y(1 - x)}}} \right]\\
&\qquad+ {\rm Li}_2 \left[ {\frac{x}{y}\frac{{(1 - y)}}{{(1 - x)}}} \right] + {\rm Li}_2 \left[ {\frac{{1 - y}}{{1 - x}}} \right] + \frac{1}{2}\log ^2 y\,.
\end{split}
\end{equation}

Putting $x=-1$ and $y=1/2$ in Kummer's formula gives

\begin{equation}\label{equ.nptw12y}
\log ^2 2 = 2\,{\rm Li}_2 \left[ { - \frac{1}{8}} \right] - 4\,{\rm Li}_2 \left[ {\frac{1}{4}} \right] - 4\,{\rm Li}_2 \left[ { - \frac{1}{2}} \right]
\end{equation}

Putting $x=\exp(i\pi/3)$ and $y=1/2$ in Kummer's formula gives 

\begin{equation}\label{equ.whrutog}
\pi ^2  = 72\,{\rm Re\; Li}_2 \left[ {\frac{1}{2}\exp \left( {\frac{{\pi i}}{3}} \right)} \right] - 18\,{\rm Li}_2 \left[ {\frac{1}{4}} \right]
\end{equation}

Putting $x=1/2$ and $y=\exp(i\pi/2)$ in Kummer's formula gives 

\begin{equation}\label{equ.c5ut706}
\begin{split}
\log ^2 2 &= 2\,{\rm Li}_2 \left[ { - \frac{1}{4}} \right] - 4\,{\rm Re\; Li}_{\rm 2} \left[ {\frac{1}{{\sqrt 2 }}\exp \left( {\frac{{3\pi i}}{4}} \right)} \right]\\
&\qquad- 4\,{\rm Re\; Li}_{\rm 2} \left[ {\frac{1}{{\left( {\sqrt 2 } \right)^3 }}\exp \left( {\frac{{\pi i}}{4}} \right)} \right]
\end{split}
\end{equation}

Putting $x=-1$ and $y=(1+i)/2$ in Kummer's formula and taking real and imaginary parts give

\begin{eqnarray}\label{equ.jkxe2jb}
&&\frac{\pi \log 2}{8} =2\,{\rm Im\,Li}_2 \left[ {\frac{1}{2}\exp \left( {\frac{{i\pi }}{2}} \right)} \right] - 2\,{\rm Im\,Li}_2 \left[ {\frac{1}{{2\sqrt 2 }}\exp \left( {\frac{{i\pi }}{4}} \right)} \right] \nonumber\\ 
&& \qquad\qquad - \,{\rm Im\,Li}_2 \left[ {\frac{1}{{4\sqrt 2 }}\exp \left( {\frac{{i\pi }}{4}} \right)} \right] 
 \end{eqnarray}
 
and

\begin{eqnarray}\label{equ.firdcin}
&&\frac{\pi^2}{32}-\frac{1}{8} \log^2 2 =2\,{\rm Re\,Li}_2 \left[ {\frac{1}{2}\exp \left( {\frac{{i\pi }}{2}} \right)} \right] +2\,{\rm Re\,Li}_2 \left[ {\frac{1}{{2\sqrt 2 }}\exp \left( {\frac{{i\pi }}{4}} \right)} \right] \nonumber\\ 
&& \qquad\qquad\qquad - \,{\rm Re\,Li}_2 \left[ {\frac{1}{{4\sqrt 2 }}\exp \left( {\frac{{i\pi }}{4}} \right)} \right]\,. 
 \end{eqnarray}

\bigskip

Another two-variable functional equation for dilogarithms, due to Abel (\mbox{Eq.~A.2.1.16} of~\cite{lewin81}) is

\begin{equation}\label{equ.u7rchcl}
\begin{split}
{\rm Li}_2 \left[ {\frac{x}{{1 - x}} \cdot \frac{y}{{1 - y}}} \right] &= {\rm Li}_2 \left[ {\frac{x}{{(1 - y)}}} \right] + {\rm Li}_2 \left[ {\frac{y}{{(1 - x)}}} \right]\\
&\qquad- {\rm Li}_2 \left[ x \right] - {\rm Li}_2 \left[ y \right] - \log (1 - x)\log (1 - y)\,.
\end{split}
\end{equation}

Putting $x=i$ and $y=-1$ in Abels formula and taking real and imaginary parts, gives

\begin{equation}\label{equ.lgahnv1}
\begin{split}
\frac{{5\pi ^2 }}{{48}} - \frac{1}{2}\log ^2 2 &= {\rm Re\; Li}_{\rm 2} \left[ {\frac{1}{{\left( {\sqrt 2 } \right)^3 }}\exp \left( {\frac{{\pi i}}{4}} \right)} \right]\\
&\qquad- {\rm Re\; Li}_{\rm 2} \left[ {\frac{1}{2}\exp \left( {\frac{{\pi i}}{2}} \right)} \right] - {\rm Re\; Li}_{\rm 2} \left[ {\frac{1}{{\sqrt 2 }}\exp \left( {\frac{{3\pi i}}{4}} \right)} \right]
\end{split}
\end{equation}

and

\begin{equation}\label{equ.in46ali}
\begin{split}
{\rm G} - \frac{{\pi \log 2}}{4} &= {\rm Im\; Li}_{\rm 2} \left[ {\frac{1}{2}\exp \left( {\frac{{\pi i}}{2}} \right)} \right]\\
&\qquad + {\rm Im\; Li}_2 \left[ {\frac{1}{{\left( {\sqrt 2 } \right)^3 }}\exp \left( {\frac{{\pi i}}{4}} \right)} \right] - {\rm Im\; Li}_2 \left[ {\frac{1}{{\sqrt 2 }}\exp \left( {\frac{{3\pi i}}{4}} \right)} \right]
\end{split}
\end{equation}

Putting $x=1/2$, $y=\exp(i\pi/3)$ in Abels formula and taking the imaginary part, gives

\begin{equation}\label{equ.ninqhk4}
5\,{\rm Cl}_2 \left( {\frac{\pi }{3}} \right) - \pi \log 2 = 6\,{\rm Im\, Li}_2 \left[ {\frac{1}{2}\exp \left( {\frac{{\pi i}}{3}} \right)} \right]
\end{equation}

Putting $x=i=y$ in Abels formula and taking real and imaginary parts, gives

\begin{equation}\label{equ.t0sjqhn}
\frac{{5\pi ^2 }}{{48}} - \frac{{\log ^2 2}}{4} = {\rm Re\; Li}_2 \left[ {\frac{1}{2}\exp \left( {\frac{{\pi i}}{2}} \right)} \right] - 2\,{\rm Re\; Li}_2 \left[ {\frac{1}{{\sqrt 2 }}\exp \left( {\frac{{3\pi i}}{4}} \right)} \right]
\end{equation}

and

\begin{equation}\label{equ.gmbyuaj}
2\,{\rm G} - \frac{{\pi \log 2}}{4} = 2\,{\rm Im\; Li}_2 \left[ {\frac{1}{{\sqrt 2 }}\exp \left( {\frac{{3\pi i}}{4}} \right)} \right] + {\rm Im\; Li}_2 \left[ {\frac{1}{2}\exp \left( {\frac{{\pi i}}{2}} \right)} \right]
\end{equation}

\subsection{Base $2^{12}$ Binary BBP-type Formulas}

Solving Eqs.~\eqref{equ.lgahnv1} and \eqref{equ.t0sjqhn} simultaneously, we find

\begin{equation}\label{equ.meiyr2w}
\begin{split}
\log ^2 2 &= 8\,{\rm Re\; Li}_2 \left[ {\frac{1}{2}\exp \left( {\frac{{\pi i}}{2}} \right)} \right] - 4\,{\rm Re\; Li}_2 \left[ {\frac{1}{{\sqrt 2 }}\exp \left( {\frac{{3\pi i}}{4}} \right)} \right]\\
&\qquad- 4\,{\rm Re\; Li}_2 \left[ {\frac{1}{{(\sqrt 2 )^3 }}\exp \left( {\frac{{\pi i}}{4}} \right)} \right]
\end{split}
\end{equation}

and

\begin{equation}\label{equ.dmg73yi}
\begin{split}
\pi^2 &= \frac{144}{5}\,{\rm Re\; Li}_2 \left[ {\frac{1}{2}\exp \left( {\frac{{\pi i}}{2}} \right)} \right] - \frac{144}{5}\,{\rm Re\; Li}_2 \left[ {\frac{1}{{\sqrt 2 }}\exp \left( {\frac{{3\pi i}}{4}} \right)} \right]\\
&\qquad- \frac{48}{5}\,{\rm Re\; Li}_2 \left[ {\frac{1}{{(\sqrt 2 )^3 }}\exp \left( {\frac{{\pi i}}{4}} \right)} \right]
\end{split}
\end{equation}

Solving Eqs.~\eqref{equ.in46ali} and \eqref{equ.gmbyuaj} simultaneously, we find

\begin{equation}\label{equ.nuiitio}
{\rm G} = 3\,{\rm Im\; Li}_2 \left[ {\frac{1}{{\sqrt 2 }}\exp \left( {\frac{{3\pi i}}{4}} \right)} \right] - {\rm Im\; Li}_2 \left[ {\frac{1}{{(\sqrt 2 )^3 }}\exp \left( {\frac{{\pi i}}{4}} \right)} \right]
\end{equation}

and

\begin{equation}\label{equ.ibwzeu0}
\begin{split}
\pi \log 2 &= 16\,{\rm Im\; Li}_2 \left[ {\frac{1}{{\sqrt 2 }}\exp \left( {\frac{{3\pi i}}{4}} \right)} \right] - 4\,{\rm  Im\; Li}_2 \left[ {\frac{1}{2}\exp \left( {\frac{{\pi i}}{2}} \right)} \right]\\
&\qquad- 8\,{\mathop{\rm Im}\nolimits} {\rm Li}_2 \left[ {\frac{1}{{(\sqrt 2 )^3 }}\exp \left( {\frac{{\pi i}}{4}} \right)} \right]
\end{split}
\end{equation}

The identities \eqref{equ.meiyr2w}, \eqref{equ.dmg73yi}, \eqref{equ.nuiitio} and \eqref{equ.ibwzeu0} facilitate the derivation of base $2^{12}$, length $24$ BBP-type formulas for the respective polylogarithm constants through the prescription of section~\ref{sec.scheme}. The explicit formulas or their variants are listed in the BBP Compendium.

\subsection{Base $2^{60}$ length $120$ Formulas}

Solving Eqs.~\eqref{equ.l71m5vv} and \eqref{equ.firdcin} simultaneously, we have

\begin{eqnarray}\label{equ.nnibiu2}
 \pi ^2  &=& 192\,{\rm Re\,Li}_2 \left[ {\frac{1}{2}\exp \left( {\frac{{i\pi }}{2}} \right)} \right] + 192\,{\rm Re\,Li}_2 \left[ {\frac{1}{{2\sqrt 2 }}\exp \left( {\frac{{i\pi }}{4}} \right)} \right] \nonumber\\ 
&&  - 96\, {\rm Re\,Li}_2 \left[ {\frac{1}{{4\sqrt 2 }}\exp \left( {\frac{{i\pi }}{4}} \right)} \right] - 24\,{\rm Li}_2 \left[ {\frac{1}{2}} \right] 
 \end{eqnarray}

and

\begin{eqnarray}\label{equ.f1143te}
 \log ^2 2 &=& 32\,{\rm Re\,Li}_2 \left[ {\frac{1}{2}\exp \left( {\frac{{i\pi }}{2}} \right)} \right] + 32\,{\rm Re\,Li}_2 \left[ {\frac{1}{{2\sqrt 2 }}\exp \left( {\frac{{i\pi }}{4}} \right)} \right] \nonumber\\ 
&&  - 16\, {\rm Re\,Li}_2 \left[ {\frac{1}{{4\sqrt 2 }}\exp \left( {\frac{{i\pi }}{4}} \right)} \right] - 6\,{\rm Li}_2 \left[ {\frac{1}{2}} \right]\,.
 \end{eqnarray}

Using Eq.~\eqref{equ.jkxe2jb} in Eq.~\eqref{equ.in46ali}, we have

\begin{eqnarray}\label{equ.eteru01}
G &=& 5\,{\rm Im\,Li}_2 \left[ {\frac{1}{2}\exp \left( {\frac{{i\pi }}{2}} \right)} \right] - {\rm Im\,Li}_2 \left[ {\frac{1}{{\sqrt 2 }}\exp \left( {\frac{{i3\pi }}{4}} \right)} \right] \nonumber\\ 
&& - 3\,{\rm Im\,Li}_2 \left[ {\frac{1}{{2\sqrt 2 }}\exp \left( {\frac{{i\pi }}{4}} \right)} \right] - 2\,{\rm Im\,Li}_2 \left[ {\frac{1}{{4\sqrt 2 }}\exp \left( {\frac{{i\pi }}{4}} \right)} \right] 
\end{eqnarray}

Applying the prescriptions of section~\ref{sec.scheme} to Eqs.~\eqref{equ.nnibiu2}, \eqref{equ.f1143te}, \eqref{equ.jkxe2jb} and \eqref{equ.eteru01}, respectively, we obtain the following base $2^{60}$, length $120$ binary BBP-type formulas:

\begin{eqnarray}
\pi^2&=&\frac{3}{2^{54}}P(2,2^{60},120,(0,-2^{58},3^2\cdot2^{58},-3^2\cdot2^{57},-5^2\cdot2^{56},-2^{56},0,\nonumber\\
&& 7\cdot2^{55},-3^2\cdot2^{55},-2^{54},0,-3^3\cdot2^{53},0,-2^{52},7\cdot2^{51},7\cdot2^{51},0,-2^{50},0,2^{53},\nonumber\\
&& 3^2\cdot2^{49},-2^{48},0,5^2\cdot2^{47},5^2\cdot2^{46},-2^{46},3^2\cdot2^{46},-3^2\cdot2^{45},0,-2^{44},0,7\cdot2^{43},\nonumber\\
&& -3^2\cdot2^{43},-2^{42},-5^2\cdot2^{41},-3^3\cdot2^{41},0,-2^{40},-3^2\cdot2^{40},-3^2\cdot2^{40},0,-2^{38},0,\nonumber\\
&& -3^2\cdot2^{37},-7\cdot2^{36},-2^{36},0,5^2\cdot2^{35},0,-2^{34},3^2\cdot2^{34},-3^2\cdot2^{33},0,-2^{32},5^2\cdot2^{31},\nonumber\\
&& 7\cdot2^{31},-3^2\cdot2^{31},-2^{30},0,-2^{30},0,-2^{28},-3^2\cdot2^{28},7\cdot2^{27},5^2\cdot2^{26},-2^{26},0,\nonumber\\
&& -3^2\cdot2^{25},3^2\cdot2^{25},-2^{24},0,5^2\cdot2^{23},0,-2^{22},-7\cdot2^{21},-3^2\cdot2^{21},0,-2^{20},0,\nonumber\\
&& -3^2\cdot2^{20},-3^2\cdot2^{19},-2^{18},0,-3^3\cdot2^{17},-5^2\cdot2^{16},-2^{16},-3^2\cdot2^{16},7\cdot2^{15},0,\nonumber\\
&& -2^{14},0,-3^2\cdot2^{13},3^2\cdot2^{13},-2^{12},5^2\cdot2^{11},5^2\cdot2^{11},0,-2^{10},3^2\cdot2^{10},2^{13},\nonumber\\
&& 0,-2^{8},0,7\cdot2^{7},7\cdot2^{6},-2^{6},0,-3^3\cdot2^{5},0,-2^{4},-3^2\cdot2^{4},7\cdot2^{3},0,-2^{2},\nonumber\\
&& -5^2\cdot2,-3^2\cdot2,3^2\cdot2,-1,0,0))\,,
\end{eqnarray}

\begin{eqnarray}
&& \log^2 2 =\nonumber\\
&& \frac{1}{2^{57}}P(2,2^{60},120,(0,-3 \cdot 2^{59},3^2 \cdot 2^{60},-19 \cdot 2^{58},-5^2 \cdot 2^{58},-3 \cdot 2^{57},0,\nonumber\\
&& 13 \cdot 2^{56},-3^2 \cdot 2^{57},-3 \cdot 2^{55},0,-5 \cdot 11 \cdot 2^{54},0,-3 \cdot 2^{53},7 \cdot 2^{53},13 \cdot 2^{52},0,-3 \cdot 2^{51},0,\nonumber\\
&& 31 \cdot 2^{50},3^2 \cdot 2^{51},-3 \cdot 2^{49},0,7^2 \cdot 2^{48},5^2 \cdot 2^{48},-3 \cdot 2^{47},3^2 \cdot 2^{48},-19 \cdot 2^{46},0,\nonumber\\
&& -3 \cdot 2^{45},0,13 \cdot 2^{44},-3^2 \cdot 2^{45},-3 \cdot 2^{43},-5^2 \cdot 2^{43},-5 \cdot 11 \cdot 2^{42},0,-3 \cdot 2^{41},-3^2 \cdot 2^{42},\nonumber\\
&& -37 \cdot 2^{40},0,-3 \cdot 2^{39},0,-19 \cdot 2^{38},-7 \cdot 2^{38},-3 \cdot 2^{37},0,7^2 \cdot 2^{36},0,-3 \cdot 2^{35},3^2 \cdot 2^{36},\nonumber\\
&& -19 \cdot 2^{34},0,-3 \cdot 2^{33},5^2 \cdot 2^{33},13 \cdot 2^{32},-3^2 \cdot 2^{33},-3 \cdot 2^{31},0,-5 \cdot 2^{30},0,-3 \cdot 2^{29},\nonumber\\
&& -3^2 \cdot 2^{30},13 \cdot 2^{28},5^2 \cdot 2^{28},-3 \cdot 2^{27},0,-19 \cdot 2^{26},3^2 \cdot 2^{27},-3 \cdot 2^{25},0,7^2 \cdot 2^{24},0,\nonumber\\
&& -3 \cdot 2^{23},-7 \cdot 2^{23},-19 \cdot 2^{22},0,-3 \cdot 2^{21},0,-37 \cdot 2^{20},-3^2 \cdot 2^{21},-3 \cdot 2^{19},0,-5 \cdot 11 \cdot 2^{18},\nonumber\\
&& -5^2 \cdot 2^{18},-3 \cdot 2^{17},-3^2 \cdot 2^{18},13 \cdot 2^{16},0,-3 \cdot 2^{15},0,-19 \cdot 2^{14},3^2 \cdot 2^{15},-3 \cdot 2^{13},\nonumber\\
&& 5^2 \cdot 2^{13},7^2 \cdot 2^{12},0,-3 \cdot 2^{11},3^2 \cdot 2^{12},31 \cdot 2^{10},0,-3 \cdot 2^{9},0,13 \cdot 2^{8},7 \cdot 2^{8},-3 \cdot 2^{7},\nonumber\\
&& 0,-5 \cdot 11 \cdot 2^{6},0,-3 \cdot 2^{5},-3^2 \cdot 2^{6},13 \cdot 2^{4},0,-3 \cdot 2^{3},-5^2 \cdot 2^{3},-19 \cdot 2^{2},\nonumber\\
&& 3^2 \cdot 2^{3},-3 \cdot 2,0,-1))\,,
\end {eqnarray}

\begin{eqnarray}
&& \pi\log 2=\nonumber\\
&&\frac{1}{2^{55}}P(2,2^{60},120,(0,2^{60},-3^2\cdot2^{57},0,-5^2\cdot2^{55},-13\cdot2^{56},0,0,-3^2\cdot2^{54},-17\cdot2^{53},\nonumber\\
&& 0,0,0,-2^{54},-7\cdot2^{50},0,0,13\cdot2^{50},0,0,3^2\cdot2^{48},-2^{50},0,0,5^2\cdot2^{45},2^{48},-3^2\cdot2^{45},0,0,\nonumber\\
&& -2^{43},0,0,-3^2\cdot2^{42},2^{44},5^2\cdot2^{40},0,0,-2^{42},3^2\cdot2^{39},0,0,13\cdot2^{38},0,0,-7\cdot2^{35},-2^{38},\nonumber\\
&& 0,0,0,-17\cdot2^{33},-3^2\cdot2^{33},0,0,-13\cdot2^{32},-5^2\cdot2^{30},0,-3^2\cdot2^{30},2^{32},0,0,0,-2^{30},\nonumber\\
&& 3^2\cdot2^{27},0,5^2\cdot2^{25},13\cdot2^{26},0,0,3^2\cdot2^{24},17\cdot2^{23},0,0,0,2^{24},7\cdot2^{20},0,0,-13\cdot2^{20},\nonumber\\
&& 0,0,-3^2\cdot2^{18},2^{20},0,0,-5^2\cdot2^{15},-2^{18},3^2\cdot2^{15},0,0,2^{13},0,0,3^2\cdot2^{12},-2^{14},-5^2\cdot2^{10},\nonumber\\
&& 0,0,2^{12},-3^2\cdot2^{9},0,0,-13\cdot2^{8},0,0,7\cdot2^{5},2^{8},0,0,0,17\cdot2^{3},3^2\cdot2^{3},0,0,\nonumber\\
&& 13\cdot2^{2},5^{2},0,3^{2},-2^{2},0,0))
\end{eqnarray}

and

\begin{eqnarray}
&& G=\nonumber\\
&& \frac{1}{2^{60}}P(2,2^{60},120,(-2^{59},3\cdot7\cdot2^{59},-7\cdot2^{60},0,-7^2\cdot2^{57},-3\cdot2^{61},2^{56},0,\nonumber\\
&& -7\cdot2^{57},-29\cdot2^{55},-2^{54},0,2^{53},-3\cdot7\cdot2^{53},-11\cdot2^{53},0,-2^{51},3\cdot2^{55},-2^{50},0,\nonumber\\
&& 7\cdot2^{51},-3\cdot7\cdot2^{49},2^{48},0,7^2\cdot2^{47},3\cdot7\cdot2^{47},-7\cdot2^{48},0,2^{45},2^{46},2^{44},0,\nonumber\\
&& -7\cdot2^{45},3\cdot7\cdot2^{43},7^2\cdot2^{42},0,2^{41},-3\cdot7\cdot2^{41},7\cdot2^{42},0,-2^{39},3\cdot2^{43},-2^{38},\nonumber\\
&& 0,-11\cdot2^{38},-3\cdot7\cdot2^{37},2^{36},0,-2^{35},-29\cdot2^{35},-7\cdot2^{36},0,2^{33},-3\cdot2^{37},-7^2\cdot2^{32},\nonumber\\
&& 0,-7\cdot2^{33},3\cdot7\cdot2^{31},-2^{30},0,2^{29},-3\cdot7\cdot2^{29},7\cdot2^{30},0,7^2\cdot2^{27},3\cdot2^{31},-2^{26},0,\nonumber\\
&& 7\cdot2^{27},29\cdot2^{25},2^{24},0,-2^{23},3\cdot7\cdot2^{23},11\cdot2^{23},0,2^{21},-3\cdot2^{25},2^{20},0,-7\cdot2^{21},\nonumber\\
&& 3\cdot7\cdot2^{19},-2^{18},0,-7^2\cdot2^{17},-3\cdot7\cdot2^{17},7\cdot2^{18},0,-2^{15},-2^{16},-2^{14},0,7\cdot2^{15},\nonumber\\
&& -3\cdot7\cdot2^{13},-7^2\cdot2^{12},0,-2^{11},3\cdot7\cdot2^{11},-7\cdot2^{12},0,2^{9},-3\cdot2^{13},2^{8},0,11\cdot2^{8},\nonumber\\
&& 3\cdot7\cdot2^{7},-2^{6},0,2^{5},29\cdot2^{5},7\cdot2^{6},0,-2^{3},3\cdot2^{7},7^2\cdot2^{2},0,\nonumber\\
&& 7\cdot2^{3},-3\cdot7\cdot2,1,0))\,.
\end{eqnarray}

\section{Degree $3$ Formulas}

No {\em proved} explicit Digit Extraction BBP-type formulas are known for $\pi^3$ and $\pi\log^2 2$. In what follows, we now present, together with their proofs, new binary degree $3$ BBP-type formulas for these and the remaining three trilogarithm constants.

\subsection{Generators of Degree $3$ BBP-type formulas}

A functional equation for trilogarithms (\mbox{Eq. A.2.6.10} of~\cite{lewin81}) reads

\begin{equation}\label{equ.subbefj}
\begin{split}
{\rm Li}_3 \left[ {\frac{{1 - x}}{{1 + x}}} \right]& - {\rm Li}_3 \left[ {\frac{{x - 1}}{{x + 1}}} \right] = 2\,{\rm Li}_3 \left[ {1 - x} \right] + 2\,{\rm Li}_3 \left[ {\frac{1}{{1 + x}}} \right]\\
&\qquad\qquad\qquad\qquad- \frac{1}{2}\,{\rm Li}_3 \left[ {1 - x^2 } \right] - \frac{7}{4}\,\zeta (3) \\
&\qquad\qquad\qquad\qquad\qquad+\frac{{\pi ^2 }}{6}\log (1 + x) - \frac{1}{3}\log ^3 (1 + x)\,.
\end{split}
\end{equation}

The use of $x=1$ in the functional equation~\eqref{equ.subbefj} gives the well-known formula

\begin{equation}\label{equ.iex39kd}
\frac{7}{8}\zeta (3) - \frac{{\pi ^2 \log 2}}{{12}} + \frac{{\log ^3 2}}{6} = {\rm Li}_3 \left[ {\frac{1}{2}} \right]\,.
\end{equation}

Another functional identity for trilogarithms (\mbox{Eq. A.2.6.11} of~\cite{lewin81}) is

\begin{equation}\label{equ.r3gvxp0}
\begin{split}
&{\rm Li}_3 \left[ {\frac{{x(1 - y)^2 }}{{y(1 - x)^2 }}} \right] + {\rm Li}_3 \left[ {xy} \right] + {\rm Li}_3 \left[ {\frac{x}{y}} \right]\\
&- 2{\rm Li}_3 \left[ {\frac{{x(1 - y)}}{{y(1 - x)}}} \right] - 2{\rm Li}_3 \left[ {\frac{{x(1 - y)}}{{(x - 1)}}} \right] - 2{\rm Li}_3 \left[ {\frac{{1 - y}}{{1 - x}}} \right]\\
& - 2{\rm Li}_3 \left[ {\frac{{(1 - y)}}{{y(x - 1)}}} \right] - 2{\rm Li}_3 \left[ x \right] - 2{\rm Li}_3 \left[ y \right] + 2\zeta (3)\\
&\qquad= \log ^2 y\log \left( {\frac{{1 - y}}{{1 - x}}} \right) - \frac{1}{3}\pi ^2 \log y - \frac{1}{3}\log ^3 y\,.
\end{split}
\end{equation}

The use of $x=-1$, $y=i$ in the above equation gives

\begin{equation}\label{equ.reeuf8k}
\frac{{35}}{8}\zeta (3) - \frac{5\pi^2\log 2}{24}  + \frac{1}{6}\log ^3 2 = 8\,{\rm Re\; Li}_3 \left[ {\frac{1}{\sqrt 2}\exp \left(\frac{i\pi}{4}\right)} \right]\,.
\end{equation}

Plugging $x=-i$, $y=1-i$ in Eq.~\eqref{equ.r3gvxp0} and taking real and imaginary parts gives

\begin{equation}\label{equ.hsmmff3}
\begin{split}
&\frac{7}{{16}}\zeta (3) + \frac{{5\pi ^2 \log 2}}{{192}} - \frac{{7\log ^3 2}}{{48}}\\
&\qquad= 2\,{\rm Re\;Li}_3 \left[ {\frac{1}{{\sqrt 2 }}\exp \left( {\frac{{i\pi }}{4}} \right)} \right] + \frac{9}{4}\,{\rm Re\;Li}_3 \left[ {\frac{1}{2}\exp \left( {\frac{{i\pi }}{2}} \right)} \right]\\
&\qquad\qquad- {\rm Re\;Li}_3 \left[ {\frac{1}{{2\sqrt 2 }}\exp \left( {\frac{{i\pi }}{4}} \right)} \right]
\end{split}
\end{equation}

and

\begin{equation}\label{equ.hv0t3gq}
\begin{split}
\frac{{3\pi \log ^2 2}}{{32}} - \frac{{9\pi ^3 }}{{128}} &= {\rm Im\; Li}_3 \left[ {\frac{1}{{2\sqrt 2 }}\exp \left( {\frac{{i\pi }}{4}} \right)} \right]\\
&\qquad+ \frac{9}{4}\,{\rm Im\; Li}_3 \left[ {\frac{1}{2}\exp \left( {\frac{{i\pi }}{2}} \right)} \right] - 6\,{\rm Im\; Li}_3 \left[ {\frac{1}{{\sqrt 2 }}\exp \left( {\frac{{i\pi }}{4}} \right)} \right]\,.
\end{split}
\end{equation}

Yet another functional equation for trilogarithms (Eq.~6.96 p 174 of~\cite{lewin81}) is 

\begin{equation}\label{equ.ww0fu0t}
\begin{split}
&{\rm Li}_3 \left[ {\frac{{x(1 - y)^2 }}{{y(1 - x)^2 }}} \right]{\rm  = 2Li}_3 \left[ {\frac{{ - x(1 - y)}}{{(1 - x)}}} \right] + 2{\rm Li}_3 \left[ {\frac{{x(1 - y)}}{{y(1 - x)}}} \right]\\
&\quad + 2{\rm Li}_3 \left[ {\frac{{ - y(1 - x)}}{{(1 - y)}}} \right] + 2{\rm Li}_3 \left[ {\frac{{1 - x}}{{1 - y}}} \right] + 2{\rm Li}_3 \left[ x \right]\\
&\qquad + 2{\rm Li}_3 \left[ y \right] - 2{\rm Li}_3 \left[ {xy} \right] - {\rm Li}_3 \left[ {\frac{x}{y}} \right] - 2\zeta (3)\\
&\qquad\quad + \frac{{\pi ^2 }}{3}\log \left( {\frac{{1 - y}}{{1 - x}}} \right) + \log ( - y)\log ^2 (1 - y)\\
&\qquad\qquad - \log ^2 ( - y)\log (1 - x) - \frac{1}{3}\log ^3 \left( {\frac{{1 - y}}{{1 - x}}} \right)\\
&\qquad\qquad\quad + \frac{1}{3}\log ^3 \left( {\frac{{ - y(1 - x)}}{{1 - y}}} \right) + \frac{1}{3}\log ^3 \left( { - \frac{{1 - y}}{y}} \right)\\
& \qquad\qquad\quad\qquad - \frac{1}{3}\log ^3 (1 - y)\,.
\end{split}
\end{equation}

Using $x=-1$, $y=1+i$ in Eq.~\eqref{equ.ww0fu0t}, simplifying and taking real and imaginary parts gives

\begin{equation}\label{equ.w2k1gjt}
\begin{split}
&\frac{7}{2}\,\zeta (3) - \frac{{15\pi ^2 \log 2}}{{64}} + \frac{{5\log ^3 2}}{{16}}\\
&\quad= 4\,{\rm Re\;Li}_3 \left[ {\frac{1}{{\sqrt 2 }}\exp \left( {\frac{{i\pi }}{4}} \right)} \right] + 4\,{\rm Re\;Li}_3 \left[ {\frac{1}{{2\sqrt 2 }}\exp \left( {\frac{{i\pi }}{4}} \right)} \right]\\
&\qquad+ \frac{7}{2}\,{\rm Re\; Li}_3 \left[ {\frac{1}{2}\exp \left( {\frac{{i\pi }}{2}} \right)} \right] - {\rm Re\;Li}_3 \left[ {\frac{1}{{4\sqrt 2 }}\exp \left( {\frac{{i\pi }}{4}} \right)} \right]
\end{split}
\end{equation}

and

\begin{equation}\label{equ.g3224vn}
\begin{split}
&\frac{{13\pi ^3 }}{{128}} - \frac{{7\pi \log ^2 2}}{{32}} \\
& \quad =4\,{\rm Im\; Li}_3 \left[ {\frac{1}{{\sqrt 2 }}\exp \left( {\frac{{i\pi }}{4}} \right)} \right] - 4\,{\rm Im\; Li}_3 \left[ {\frac{1}{{2\sqrt 2 }}\exp \left( {\frac{{i\pi }}{4}} \right)} \right]\\
& \qquad+ \frac{7}{2}\,{\rm Im\; Li}_3 \left[ {\frac{1}{2}\exp \left( {\frac{{i\pi }}{2}} \right)} \right] - {\rm Im\; Li}_3 \left[ {\frac{1}{{4\sqrt 2 }}\exp \left( {\frac{{i\pi }}{4}} \right)} \right]\,.
\end{split}
\end{equation}

Using $x=-1$, $y=1/2$ in Eq.~\eqref{equ.ww0fu0t} gives the identity

\begin{equation}\label{equ.mqgomtf}
7\,\zeta (3) - \pi ^2 \log 2 + 3\,\log ^3 2 = 9\,{\rm L}_3 \left[ {\frac{1}{4}} \right] - 2\,{\rm L}_3 \left[ { - \frac{1}{8}} \right]\,.
\end{equation}

We are now ready to derive formulas for $\pi^3$, $\pi\log^2 2$, $\zeta(3)$, $\pi^2\log 2$ and $\log^3 2$.

\subsection{Base $2^{12}$ Formulas}

Solving Eqs.~\eqref{equ.iex39kd}, \eqref{equ.reeuf8k} and \eqref{equ.mqgomtf} simultaneously, we find the BBP-ready identities:

\begin{eqnarray}\label{equ.ls90hkt}
 \zeta (3) &=& \frac{{128}}{{21}}\,{\rm Re\,Li}_3 \left[ {\frac{1}{{\sqrt 2 }}\exp \left( {\frac{{i\pi }}{4}} \right)} \right] - \frac{{88}}{{21}}\,{\rm Li}_3 \left[ {\frac{1}{2}} \right] + \frac{{12}}{7}\,{\rm Li}_3 \left[ {\frac{1}{4}} \right]\nonumber \\ 
&& \qquad - \frac{8}{{21}}\,{\rm Li}_3 \left[ { - \frac{1}{8}} \right] 
\end{eqnarray}

\begin{eqnarray}\label{equ.vjkkvn0}
\pi^2 \log 2 &=& \frac{{320}}{{3}}\,{\rm Re\,Li}_3 \left[ {\frac{1}{{\sqrt 2 }}\exp \left( {\frac{{i\pi }}{4}} \right)} \right] - \frac{{328}}{{3}}\,{\rm Li}_3 \left[ {\frac{1}{2}} \right] + 48\,{\rm Li}_3 \left[ {\frac{1}{4}} \right]\nonumber \\ 
&& \qquad - \frac{32}{{3}}\,{\rm Li}_3 \left[ { - \frac{1}{8}} \right] 
\end{eqnarray}

and

\begin{eqnarray}\label{equ.pr6an83}
 \log^3 2 &=& \frac{{64}}{{3}}\,{\rm Re\,Li}_3 \left[ {\frac{1}{{\sqrt 2 }}\exp \left( {\frac{{i\pi }}{4}} \right)} \right] - \frac{{80}}{{3}}\,{\rm Li}_3 \left[ {\frac{1}{2}} \right] + 15\,{\rm Li}_3 \left[ {\frac{1}{4}} \right]\nonumber \\ 
&& \qquad - \frac{10}{{3}}\,{\rm Li}_3 \left[ { - \frac{1}{8}} \right]\,. 
\end{eqnarray}

Application of the prescriptions of section~\ref{sec.scheme} to the identities \eqref{equ.ls90hkt}, \eqref{equ.vjkkvn0} and \eqref{equ.pr6an83} facilitates the derivation of base $2^{12}$, length $24$ BBP-type formulas for the respective polylogarithm constants. The explicit formulas or their variants are as listed in the BBP Compendium.

\subsection{Base $2^{60}$ Formulas}

Eliminating $\pi\log^2 2$ between Eqs.~\eqref{equ.hv0t3gq} and \eqref{equ.g3224vn}, we have the following BBP-ready formula for $\pi^3$:

\begin{equation}\label{equ.qfqfmla}
\begin{split}
\pi ^3  &= 16\,{\rm Im\; Li}_3 \left[ {\frac{1}{{4\sqrt 2 }}\exp \left( {\frac{{i\pi }}{4}} \right)} \right] + 160\,{\rm Im\; Li}_3 \left[ {\frac{1}{{\sqrt 2 }}\exp \left( {\frac{{i\pi }}{4}} \right)} \right]\\
&\qquad {\rm  + }\frac{{80}}{3}\,{\rm Im\; Li}_3 \left[ {\frac{1}{{2\sqrt 2 }}\exp \left( {\frac{{i\pi }}{4}} \right)} \right] - 140\,{\rm Im\; Li}_3 \left[ {\frac{1}{2}\exp \left( {\frac{{i\pi }}{2}} \right)} \right]\,.
\end{split}
\end{equation}

Eliminating $\pi^3$ between Eqs.~\eqref{equ.hv0t3gq} and \eqref{equ.g3224vn}, we have the following BBP-ready formula for $\pi\log^2 2$:

\begin{equation}\label{equ.mtrgxwl}
\begin{split}
\pi\log^2 2  &= 12\,{\rm Im\; Li}_3 \left[ {\frac{1}{{4\sqrt 2 }}\exp \left( {\frac{{i\pi }}{4}} \right)} \right] + 56\,{\rm Im\; Li}_3 \left[ {\frac{1}{{\sqrt 2 }}\exp \left( {\frac{{i\pi }}{4}} \right)} \right]\\
&\qquad {\rm  + }\frac{{92}}{3}\,{\rm Im\; Li}_3 \left[ {\frac{1}{{2\sqrt 2 }}\exp \left( {\frac{{i\pi }}{4}} \right)} \right] - 81\,{\rm Im\; Li}_3 \left[ {\frac{1}{2}\exp \left( {\frac{{i\pi }}{2}} \right)} \right]\,.
\end{split}
\end{equation}

Eliminating $\pi^2\log 2$ and $\log^3 2$ between Eqs.~\eqref{equ.reeuf8k}, \eqref{equ.hsmmff3} and \eqref{equ.w2k1gjt} we find

\begin{equation}\label{equ.fgc5kj5}
\begin{split}
\zeta(3)  &= \frac{16}{7}\,{\rm Re\; Li}_3 \left[ {\frac{1}{{4\sqrt 2 }}\exp \left( {\frac{{i\pi }}{4}} \right)} \right] + \frac{32}{7}\,{\rm Re\; Li}_3 \left[ {\frac{1}{{\sqrt 2 }}\exp \left( {\frac{{i\pi }}{4}} \right)} \right]\\
&\qquad -\frac{{48}}{7}\,{\rm Re\; Li}_3 \left[ {\frac{1}{{2\sqrt 2 }}\exp \left( {\frac{{i\pi }}{4}} \right)} \right] - \frac{92}{7}\,{\rm Re\; Li}_3 \left[ {\frac{1}{2}\exp \left( {\frac{{i\pi }}{2}} \right)} \right]\,.
\end{split}
\end{equation}

Eliminating $\zeta(3)$ and $\log^3 2$ between Eqs.~\eqref{equ.reeuf8k}, \eqref{equ.hsmmff3} and \eqref{equ.w2k1gjt} we find

\begin{equation}\label{equ.i6p2gpx}
\begin{split}
\pi^2\log 2  &= \frac{312}{5}\,{\rm Re\; Li}_3 \left[ {\frac{1}{{4\sqrt 2 }}\exp \left( {\frac{{i\pi }}{4}} \right)} \right] + \frac{336}{5}\,{\rm Re\; Li}_3 \left[ {\frac{1}{{\sqrt 2 }}\exp \left( {\frac{{i\pi }}{4}} \right)} \right]\\
&\qquad -\frac{{904}}{5}\,{\rm Re\; Li}_3 \left[ {\frac{1}{{2\sqrt 2 }}\exp \left( {\frac{{i\pi }}{4}} \right)} \right] - \frac{1866}{5}\,{\rm Re\; Li}_3 \left[ {\frac{1}{2}\exp \left( {\frac{{i\pi }}{2}} \right)} \right]\,.
\end{split}
\end{equation}

Elimination of $\pi^2\log 2$ and $\zeta(3)$ between Eqs.~\eqref{equ.reeuf8k}, \eqref{equ.hsmmff3} and \eqref{equ.w2k1gjt} gives

\begin{equation}\label{equ.uyefz2x}
\begin{split}
\log^3 2  &= 18\,{\rm Re\; Li}_3 \left[ {\frac{1}{{4\sqrt 2 }}\exp \left( {\frac{{i\pi }}{4}} \right)} \right] + 12\,{\rm Re\; Li}_3 \left[ {\frac{1}{{\sqrt 2 }}\exp \left( {\frac{{i\pi }}{4}} \right)} \right]\\
&\qquad -46\,{\rm Re\; Li}_3 \left[ {\frac{1}{{2\sqrt 2 }}\exp \left( {\frac{{i\pi }}{4}} \right)} \right] - \frac{243}{2}\,{\rm Re\; Li}_3 \left[ {\frac{1}{2}\exp \left( {\frac{{i\pi }}{2}} \right)} \right]\,.
\end{split}
\end{equation}

As in the previous cases, the application of the prescriptions of section~\ref{sec.scheme} to the identities \eqref{equ.qfqfmla}, \eqref{equ.mtrgxwl}, \eqref{equ.fgc5kj5}, \eqref{equ.i6p2gpx} and \eqref{equ.uyefz2x} facilitates the derivation of base $2^{60}$, length $120$ BBP-type formulas for the respective polylogarithm constants. The explicit formulas or their variants are as listed in the BBP Compendium.

\section{Degree $4$ Formulas}

Next we derive binary BBP-type formulas for $\pi ^4$, $\pi ^2 \log ^2 2$, $\log ^4 2$ and two linear combinations of ${\rm Cl}_4(\pi/2)$, $\pi\log^3 2$ and $\pi^3\log 2$. We will also give formal proofs of the known but hitherto unproved formulas for these polylogarithm constants.

\subsection{\mbox{Generators of Degree 4 Binary BBP-type Formulas}}

The two-variable degree 4 polylogarithm functional equation (Eq.~7.90 pg. 211 of Lewin's book~\cite{lewin81}) reads 	

\begin{equation}\label{equ.ieibdj8}
\begin{split}
&\operatorname{Li}_4 \left[-{\frac {{x}^{2}y\eta}{\xi}} \right] +\operatorname{Li}_4 \left[-{\frac {{
y}^{2}x\xi}{\eta}} \right] +\operatorname{Li}_4 \left[{\frac {{x}^{2}y}{{\eta}^{2}\xi
}} \right] +\operatorname{Li}_4 \left[{\frac {{y}^{2}x}{{\xi}^{2}\eta}} \right]\\ 
&=6\,\operatorname{Li}_4\left[xy \right] +6\,\operatorname{Li}_4 \left[{\frac {xy}{\eta\,\xi}} \right] +6
\,\operatorname{Li}_4 \left[-{\frac {xy}{\eta}} \right] +6\,\operatorname{Li}_4 \left[-{\frac {xy}{
\xi}} \right]\\
&+3\,\operatorname{Li}_4 \left[x\eta \right] +3\,\operatorname{Li}_4 \left[y\xi
 \right] +3\,\operatorname{Li}_4 \left[{\frac {x}{\eta}} \right] +3\,\operatorname{Li}_4 \left[{
\frac {y}{\xi}} \right] +3\,\operatorname{Li}_4 \left[-{\frac {x\eta}{\xi}} \right]\\
& +3\,\operatorname{Li}_4 \left[-{\frac {y\xi}{\eta}} \right] +3\,\operatorname{Li}_4 \left[-{\frac {x}
{\eta\,\xi}} \right] +3\,\operatorname{Li}_4 \left[-{\frac {y}{\eta\,\xi}} \right] -6
\,\operatorname{Li}_4 \left[x \right]\\
& -6\,\operatorname{Li}_4 \left[y \right] -6\,\operatorname{Li}_4 \left[-{
\frac {x}{\xi}} \right] -6\,\operatorname{Li}_4 \left[-{\frac {y}{\eta}} \right] +3/2
\,\log^2\xi\log^2\eta\,,
\end{split}
\end{equation}

where $\xi=1-x$, $\eta=1-y$.

Evaluating Eq.~\eqref{equ.ieibdj8} at coordinates $(1/2,1/2)$ gives

\begin{equation}\label{equ.dch08yt}
\begin{split}
 &\frac{{\pi ^4 }}{9} - \frac{1}{2}\,\pi ^2 \log ^2 2 + \frac{5}{4}\,\log ^4 2 \\ 
  &\qquad= 24\,{\rm Li}_4 \left[ {\frac{1}{2}} \right] + 2\,{\rm Li}_4 \left[ { - \frac{1}{8}} \right] - \frac{{27}}{2}\,{\rm Li}_4 \left[ {\frac{1}{4}} \right]\,.
 \end{split}
\end{equation}

Evaluating the functional equation Eq.~\eqref{equ.ieibdj8} at coordinates $(i,i)$, simplifying and taking real and imaginary parts, we obtain

\begin{equation}\label{equ.cx196an}
\begin{split}
 &\frac{{349\pi ^4 }}{{9216}} - \frac{{7\pi ^2 \log ^2 2}}{{128}} + \frac{5}{{64}}\,\log ^4 2 \\ 
  &\quad= 2\,{\mathop{\rm Re}\nolimits}\, {\rm Li}_4 \left[ {\frac{1}{{2\sqrt 2 }}\exp \left( {\frac{{i\pi }}{4}} \right)} \right] + 12\,{\mathop{\rm Re}\nolimits}\, {\rm Li}_4 \left[ {\frac{1}{{\sqrt 2 }}\exp \left( {\frac{{i\pi }}{4}} \right)} \right] \\ 
  &\qquad- \frac{{27}}{4}\,{\mathop{\rm Re}\nolimits}\, {\rm Li}_4 \left[ {\frac{1}{2}\exp \left( {\frac{{i\pi }}{2}} \right)} \right] - 6\,{\rm Li}_4 \left[ {\frac{1}{2}} \right] 
\end{split}
\end{equation}

and

\begin{equation}\label{equ.px7sjmm}
\begin{split}
&20{\rm Cl}_{\rm 4} \left( {\frac{\pi }{2}} \right) + \frac{{3\pi \log ^3 2}}{{32}} - \frac{{27\pi ^3 \log 2}}{{128}} \\ 
 &\quad =  - 2\,{\mathop{\rm Im}\nolimits}\, {\rm Li}_4 \left[ {\frac{1}{{2\sqrt 2 }}\exp \left( {\frac{{i\pi }}{4}} \right)} \right] - \frac{{27}}{4}\,{\mathop{\rm Im}\nolimits} \,{\rm Li}_4 \left[ {\frac{1}{2}\exp \left( {\frac{{i\pi }}{2}} \right)} \right] \\ 
  &\qquad+ 36\,{\mathop{\rm Im}\nolimits} \,{\rm Li}_4 \left[ {\frac{1}{{\sqrt 2 }}\exp \left( {\frac{{i\pi }}{4}} \right)} \right] \,.
\end{split}
\end{equation}

Evaluating Eq.~\eqref{equ.ieibdj8} at $(1+i,1/2)$, simplifying and taking the real part yields
\begin{equation}\label{equ.dr1o6yd}
\begin{split}
 &\frac{{1697\pi ^4 }}{{9216}} - \frac{{137\pi ^2 \log ^2 2}}{{384}} + \frac{{115}}{{192}}\log ^4 2 \\ 
 &\quad =  - 14\,{\mathop{\rm Re}\nolimits}\,{\rm Li}_4 \left[ {\frac{1}{{\sqrt 2 }}\exp \left( {\frac{{3\pi i}}{4}} \right)} \right] + 30\,{\mathop{\rm Re}\nolimits}\, {\rm Li}_4 \left[ {\frac{1}{{\sqrt 2 }}\exp \left( {\frac{{i\pi }}{4}} \right)} \right] \\ 
 &\quad - 6\,{\mathop{\rm Re}\nolimits} \,{\rm Li}_4 \left[ {\frac{1}{{2\sqrt 2 }}\exp \left( {\frac{{i\pi }}{2}} \right)} \right] - 7\,{\rm Li}_4 \left[ {\frac{1}{2}} \right] + \frac{{11}}{8}\,{\rm Li}_4 \left[ { - \frac{1}{4}} \right]\,.
\end{split}
\end{equation}  

Finally, evaluating the identity at $((1-i)/2,1/2)$ and taking real and imaginary parts, we find
\begin{equation}\label{equ.ppravew}
\begin{split}
&\frac{{1265\pi ^4 }}{{18432}} - \frac{{113\pi ^2 \log ^2 2}}{{768}} + \frac{{91}}{{384}}\log ^4 2 \\ 
  &\quad= {\mathop{\rm Re}\nolimits}\, {\rm Li}_4 \left[ {\frac{1}{{4\sqrt 2 }}\exp \left( {\frac{{\pi i}}{4}} \right)} \right] - 12\,{\mathop{\rm Re}\nolimits} \,{\rm Li}_4 \left[ {\frac{1}{{2\sqrt 2 }}\exp \left( {\frac{{i\pi }}{4}} \right)} \right] \\ 
  &\qquad+ 5\,{\mathop{\rm Re}\nolimits} \,{\rm Li}_4 \left[ {\frac{1}{{\sqrt 2 }}\exp \left( {\frac{{i\pi }}{4}} \right)} \right] - 7\,{\mathop{\rm Re}\nolimits}\, {\rm Li}_4 \left[ {\frac{1}{{\sqrt 2 }}\exp \left( {\frac{{3\pi i}}{4}} \right)} \right] \\ 
 &\qquad\qquad - \frac{{95}}{8}\,{\mathop{\rm Re}\nolimits} \,{\rm Li}_4 \left[ {\frac{1}{2}\exp \left( {\frac{{\pi i}}{2}} \right)} \right] + \frac{{11}}{8}\,{\rm Li}_4 \left[ { - \frac{1}{4}} \right] + 6\,{\rm Li}_4 \left[ {\frac{1}{2}} \right]
\end{split}
\end{equation}

and

\begin{equation}\label{equ.wwh5r87}
\begin{split}
& 12\,{\rm Cl}_{\rm 4} \left( {\frac{\pi }{2}} \right) + \frac{{29\pi \log ^3 2}}{{192}} - \frac{{47\pi ^3 \log 2}}{{256}} \\ 
 &\quad=  - 6\,{\mathop{\rm Im}\nolimits} {\rm Li}_4 \left[ {\frac{1}{{2\sqrt 2 }}\exp \left( {\frac{{i\pi }}{4}} \right)} \right] + \frac{{95}}{8}\,{\mathop{\rm Im}\nolimits} \,{\rm Li}_4 \left[ {\frac{1}{2}\exp \left( {\frac{{i\pi }}{2}} \right)} \right] \\ 
 &\qquad + {\mathop{\rm Im}\nolimits} \,{\rm Li}_4 \left[ {\frac{1}{{\sqrt 2 }}\,\exp \left( {\frac{{i\pi }}{4}} \right)} \right] + 7\,{\mathop{\rm Im}\nolimits} \,{\rm Li}_4 \left[ {\frac{1}{{\sqrt 2 }}\exp \left( {\frac{{3\pi i}}{4}} \right)} \right] \\ 
  &\quad\qquad- {\mathop{\rm Im}\nolimits} \,{\rm Li}_4 \left[ {\frac{1}{{4\sqrt 2 }}\exp \left( {\frac{{i\pi }}{4}} \right)} \right]\,.
\end{split}
\end{equation}

We are now ready to derive formulas for $\pi^4$, $\pi^2\log^2 2$, $\log^4 2$, and the two linear combinations of ${\rm Cl}_4(\pi/2)$, $\pi\log^3 2$ and $\pi^3\log 2$.

\subsection{Base $2^{12}$ Binary BBP-type Formulas}

Solving Eqs.~\eqref{equ.dch08yt}, \eqref{equ.cx196an} and \eqref{equ.dr1o6yd} simultaneously, we find

\begin{equation}\label{equ.b3ywont}
\begin{split}
 \pi ^4  &= \frac{{27648}}{{41}}\,{\mathop{\rm Re}\nolimits} \,{\rm Li}_4 \left[ {\frac{1}{{2\sqrt 2 }}\exp \left( {\frac{{i\pi }}{4}} \right)} \right] + \frac{{51840}}{{41}}\,{\mathop{\rm Re}\nolimits} \,{\rm Li}_4 \left[ {\frac{1}{{\sqrt 2 }}\exp \left( {\frac{{i\pi }}{4}} \right)} \right] \\ 
  &\qquad+ \frac{{24192}}{{41}}\,{\mathop{\rm Re}\nolimits} \,{\rm Li}_4 \left[ {\frac{1}{{\sqrt 2 }}\exp \left( {\frac{{3i\pi }}{4}} \right)} \right] - \frac{{58320}}{{41}}\,{\mathop{\rm Re}\nolimits} \,{\rm Li}_4 \left[ {\frac{1}{2}\exp \left( {\frac{{i\pi }}{2}} \right)} \right] \\ 
  &\qquad- \frac{{32832}}{{41}}\,{\rm Li}_4 \left[ {\frac{1}{2}} \right] - \frac{{3888}}{{41}}\,{\rm Li}_4 \left[ {\frac{1}{4}} \right] + \frac{{576}}{{41}}\,{\rm Li}_4 \left[ { - \frac{1}{8}} \right] - \frac{{2376}}{{41}}\,{\rm Li}_4 \left[ { - \frac{1}{4}} \right] \,,
\end{split}
\end{equation}

\begin{equation}\label{equ.g0s8xgl}
\begin{split}
 \pi ^2 \log ^2 2 &= \frac{{98944}}{{123}}\,{\mathop{\rm Re}\nolimits} \,{\rm Li}_4 \left[ {\frac{1}{{2\sqrt 2 }}\exp \left( {\frac{{i\pi }}{4}} \right)} \right] + \frac{{47408}}{{41}}\,{\mathop{\rm Re}\nolimits}\, {\rm Li}_4 \left[ {\frac{1}{{\sqrt 2 }}\exp \left( {\frac{{i\pi }}{4}} \right)} \right] \\ 
  &\qquad+ \frac{{31920}}{{41}}\,{\mathop{\rm Re}\nolimits}\, {\rm Li}_4 \left[ {\frac{1}{{\sqrt 2 }}\exp \left( {\frac{{3i\pi }}{4}} \right)} \right] - \frac{{65142}}{{41}}\,{\mathop{\rm Re}\nolimits} \,{\rm Li}_4 \left[ {\frac{1}{2}\exp \left( {\frac{{i\pi }}{2}} \right)} \right] \\ 
  &\qquad- \frac{{30200}}{{41}}\,{\rm Li}_4 \left[ {\frac{1}{2}} \right] - \frac{{6606}}{{41}}\,{\rm Li}_4 \left[ {\frac{1}{4}} \right] + \frac{{2936}}{{123}}\,{\rm Li}_4 \left[ { - \frac{1}{8}} \right] - \frac{{3135}}{{41}}\,{\rm Li}_4 \left[ { - \frac{1}{4}} \right]
\end{split}
\end{equation}

and

\begin{equation}\label{equ.pezki4s}
\begin{split}
 \log ^2 4& = \frac{{161024}}{{615}}\,{\mathop{\rm Re}\nolimits} \,{\rm Li}_4 \left[ {\frac{1}{{2\sqrt 2 }}\exp \left( {\frac{{i\pi }}{4}} \right)} \right] + \frac{{71776}}{{205}}\,{\mathop{\rm Re}\nolimits}\, {\rm Li}_4 \left[ {\frac{1}{{\sqrt 2 }}\exp \left( {\frac{{i\pi }}{4}} \right)} \right] \\ 
  &\qquad+ \frac{{53088}}{{205}}\,{\mathop{\rm Re}\nolimits} \,{\rm Li}_4 \left[ {\frac{1}{{\sqrt 2 }}\exp \left( {\frac{{3i\pi }}{4}} \right)} \right] - \frac{{104364}}{{205}}\,{\mathop{\rm Re}\nolimits} \,{\rm Li}_4 \left[ {\frac{1}{2}\exp \left( {\frac{{i\pi }}{2}} \right)} \right] \\ 
  &\qquad- \frac{{41872}}{{205}}\,{\rm Li}_4 \left[ {\frac{1}{2}} \right] - \frac{{13698}}{{205}}\,{\rm Li}_4 \left[ {\frac{1}{4}} \right] + \frac{{6088}}{{615}}\,{\rm Li}_4 \left[ { - \frac{1}{8}} \right] - \frac{{5214}}{{205}}\,{\rm Li}_4 \left[ { - \frac{1}{4}} \right]\,.
\end{split}
\end{equation}

The identities \eqref{equ.px7sjmm}, \eqref{equ.b3ywont}, \eqref{equ.g0s8xgl} and \eqref{equ.pezki4s} facilitate the derivation of base $2^{12}$, length $24$ BBP-type formulas for the respective polylogarithm constants. The explicit formulas are listed in the BBP Compendium.

\subsection{Base $2^{60}$ Binary BBP-type Formulas}
Solving Eqs.~\eqref{equ.cx196an}, \eqref{equ.dr1o6yd} and \eqref{equ.ppravew} simultaneously, we find

\begin{equation}\label{equ.b2nxb0w}
\begin{split}
\pi ^4  &= \frac{{57600}}{{71}}{\mathop{\rm Re}\nolimits} {\rm Li}_4 \left[ {\frac{1}{{2\sqrt 2 }}\exp \left( {\frac{{i\pi }}{4}} \right)} \right] + \frac{{442368}}{{497}}{\mathop{\rm Re}\nolimits} {\rm Li}_4 \left[ {\frac{1}{{\sqrt 2 }}\exp \left( {\frac{{i\pi }}{4}} \right)} \right] \\ 
  &\qquad+ \frac{{13824}}{{71}}{\mathop{\rm Re}\nolimits} {\rm Li}_4 \left[ {\frac{1}{{\sqrt 2 }}\exp \left( {\frac{{3i\pi }}{4}} \right)} \right] + \frac{{239328}}{{497}}{\mathop{\rm Re}\nolimits} {\rm Li}_4 \left[ {\frac{1}{2}\exp \left( {\frac{{i\pi }}{2}} \right)} \right] \\ 
 &\qquad-\frac{{34560}}{{497}}{\mathop{\rm Re}\nolimits} {\rm Li}_4 \left[ {\frac{1}{{4\sqrt 2 }}\exp \left( {\frac{{i\pi }}{4}} \right)} \right] - \frac{{432000}}{{497}}{\rm Li}_4 \left[ {\frac{1}{2}} \right] - \frac{{4752}}{{71}}{\rm Li}_4 \left[ { - \frac{1}{4}} \right] \,,
\end{split}
\end{equation}

\begin{equation}\label{equ.nxvhffk}
\begin{split}
\pi ^2 \log ^2 2 &= \frac{{73632}}{{71}}{\mathop{\rm Re}\nolimits} {\rm Li}_4 \left[ {\frac{1}{{2\sqrt 2 }}\exp \left( {\frac{{i\pi }}{4}} \right)} \right] + \frac{{258592}}{{497}}{\mathop{\rm Re}\nolimits} {\rm Li}_4 \left[ {\frac{1}{{\sqrt 2 }}\exp \left( {\frac{{i\pi }}{4}} \right)} \right] \\ 
  &\qquad+ \frac{{7584}}{{71}}{\mathop{\rm Re}\nolimits} {\rm Li}_4 \left[ {\frac{1}{{\sqrt 2 }}\exp \left( {\frac{{3i\pi }}{4}} \right)} \right] + \frac{{818152}}{{497}}{\mathop{\rm Re}\nolimits} {\rm Li}_4 \left[ {\frac{1}{2}\exp \left( {\frac{{i\pi }}{2}} \right)} \right] \\ 
  &\qquad- \frac{{58720}}{{497}}{\mathop{\rm Re}\nolimits} {\rm Li}_4 \left[ {\frac{1}{{4\sqrt 2 }}\exp \left( {\frac{{i\pi }}{4}} \right)} \right] - \frac{{423872}}{{497}}{\rm Li}_4 \left[ {\frac{1}{2}} \right] - \frac{{6512}}{{71}}{\rm Li}_4 \left[ { - \frac{1}{4}} \right]
\end{split}
\end{equation}

and

\begin{equation}\label{equ.hs5hvdo}
\begin{split}
\log ^4 2 &= \frac{{25440}}{{71}}{\mathop{\rm Re}\nolimits} {\rm Li}_4 \left[ {\frac{1}{{2\sqrt 2 }}\exp \left( {\frac{{i\pi }}{4}} \right)} \right] + \frac{{42928}}{{497}}{\mathop{\rm Re}\nolimits} {\rm Li}_4 \left[ {\frac{1}{{\sqrt 2 }}\exp \left( {\frac{{i\pi }}{4}} \right)} \right] \\ 
  &\qquad- \frac{{1392}}{{71}}{\mathop{\rm Re}\nolimits} {\rm Li}_4 \left[ {\frac{1}{{\sqrt 2 }}\exp \left( {\frac{{3i\pi }}{4}} \right)} \right] + \frac{{413758}}{{497}}{\mathop{\rm Re}\nolimits} {\rm Li}_4 \left[ {\frac{1}{2}\exp \left( {\frac{{i\pi }}{2}} \right)} \right] \\ 
  &\qquad- \frac{{24352}}{{497}}{\mathop{\rm Re}\nolimits} {\rm Li}_4 \left[ {\frac{1}{{4\sqrt 2 }}\exp \left( {\frac{{i\pi }}{4}} \right)} \right] - \frac{{125480}}{{497}}{\rm Li}_4 \left[ {\frac{1}{2}} \right] - \frac{{2255}}{{71}}{\rm Li}_4 \left[ { - \frac{1}{4}} \right] \,.
\end{split}
\end{equation}

The identities \eqref{equ.wwh5r87}, \eqref{equ.b2nxb0w}, \eqref{equ.nxvhffk} and \eqref{equ.hs5hvdo} facilitate the derivation of base $2^{60}$, length $120$ BBP-type formulas for the respective polylogarithm constants. The explicit formulas are listed in the BBP Compendium.

\section{Degree $5$ Formulas}

Next we give formal proofs for the binary BBP-type formulas for $\zeta (5)$, $\pi ^4 \log 2$, $\pi ^2 \log ^3 2$ and $\log ^5 2$. These formulas were found experimentally by David H. Bailey, using his PSLQ algorithm. This section provides the first avenue where the hitherto unproved formulas are formally proved.

\subsection{\mbox{Generators of Degree $5$ Binary BBP-type Formulas}}

The following two-variable degree 5 polylogarithm functional equation was derived by Broadhurst~\cite{broadhurst98} (Eq.~63 pg. 5)	

\begin{equation}
\begin{split}
\label{equ.t9v474l}
&{\rm Li}_5 \left[{\frac {x\alpha}{y\beta}} \right] +{\rm Li}_5 \left[x\alpha\,y
\eta \right] +{\rm Li}_5 \left[{\frac {x\alpha\,\beta}{\eta}} \right] +{\rm Li}_5
 \left[x\xi\,y\beta \right] +{\rm Li}_5 \left[{\frac {x\xi}{y\eta}}
 \right]\\
 &+{\rm Li}_5 \left[{\frac {x\xi\,\eta}{\beta}} \right] +{\rm Li}_5 \left[{
\frac {\alpha\,y\beta}{\xi}} \right] +{\rm Li}_5 \left[{\frac {\alpha}{\xi\,
y\eta}} \right] +{\rm Li}_5 \left[{\frac {\alpha\,\eta}{\xi\,\beta}}
 \right]\\
 &-9\,{\rm Li}_5 \left[xy \right] -9\,{\rm Li}_5 \left[x\beta \right] -9\,{\rm Li}_5
 \left[x\eta \right] -9\,{\rm Li}_5 \left[{\frac {x}{y}} \right] -9\,{\rm Li}_5
 \left[{\frac {x}{\beta}} \right]\\ 
&-9\,{\rm Li}_5 \left[{\frac {x}{\eta}}
 \right] -9\,{\rm Li}_5 \left[\alpha\,y \right] -9\,{\rm Li}_5 \left[\alpha\,\beta
 \right] -9\,{\rm Li}_5 \left[\alpha\,\eta \right]\\ 
&-9\,{\rm Li}_5 \left[{\frac {
\alpha}{y}} \right] -9\,{\rm Li}_5 \left[{\frac {\alpha}{\beta}} \right] -9
\,{\rm Li}_5 \left[{\frac {\alpha}{\eta}} \right] -9\,{\rm Li}_5 \left[\xi\,y
 \right] -9\,{\rm Li}_5 \left[\xi\,\beta \right]\\
 &-9\,{\rm Li}_5 \left[\xi\,\eta
 \right] -9\,{\rm Li}_5 \left[{\frac {y}{\xi}} \right] -9\,{\rm Li}_5 \left[{
\frac {\beta}{\xi}} \right] -9\,{\rm Li}_5 \left[{\frac {\eta}{\xi}}
 \right]\\
&  +18\,{\rm Li}_5 \left[x \right] +18\,{\rm Li}_5 \left[\alpha \right] +18
\,{\rm Li}_5 \left[\xi \right] +18\,{\rm Li}_5 \left[y \right] +18\,{\rm Li}_5 \left[
\beta \right]\\
&+18\,{\rm Li}_5 \left[\eta \right] -18\,\zeta (5) =
3/10\, \left( \log\xi\right) ^{5}+3/4\, \left( \log y -\log x  \right)  \left( \log\xi  \right) ^{4}\\
&+3/2\, \left( 3\,\log y -\log\eta  \right)  \left( \log\eta\right) ^{2} \left( \log \xi  \right) ^{2}+1/2\,{\pi }
^{2} \left( \log \xi -3\,\log \eta\right)\left( \log \xi\right)^{2}+1/5\,{\pi }^{4}\log  \xi\,.
\end{split}
\end{equation}

In the above formula $\xi=1-x$, $\eta=1-y$, $\alpha=-x/\xi$ and $\beta=-y/\eta$.

\bigskip

Evaluating Eq.~\eqref{equ.t9v474l} at coordinates $(1/2,1/2)$ gives

\begin{equation}\label{equ.deg5a}
\begin{split}
 &\frac{{403}}{4}\,\zeta (5) - \frac{2}{3}\,\pi ^4 \log 2 + \pi ^2 \log ^3 2 - \frac{3}{2}\,\log ^5 2 \\ 
  &\quad = 144\,{\rm Li}_5 \left[ {\frac{1}{2}} \right] - \frac{{81}}{2}\,{\rm Li}_5 \left[ {\frac{1}{4}} \right] + 4\,{\rm Li}_5 \left[ { - \frac{1}{8}} \right]\,.
 \end{split}
\end{equation}

Evaluating at $(-i,i)$ and taking the real part gives

\begin{equation}\label{equ.deg5b}
\begin{split}
 &\frac{{4371}}{{128}}\,\zeta (5) - \frac{{349}}{{3072}}\,\pi ^4 \log 2 + \frac{7}{{128}}\,\pi ^2 \log ^3 2 - \frac{3}{{64}}\,\log ^5 2 \\ 
 &\quad = 36\,{\mathop{\rm Re}\nolimits}\, {\rm Li}_5 \left[ {\frac{1}{{\sqrt 2 }}\exp \left( {\frac{{i\pi }}{4}} \right)} \right] - 36{\mathop{\rm Re}\nolimits}\, {\rm Li}_5 \left[ {\frac{1}{{\sqrt 2 }}\exp \left( {\frac{{3i\pi }}{4}} \right)} \right] \\ 
  &\qquad+ 4{\mathop{\rm Re}\nolimits}\, {\rm Li}_5 \left[ {\frac{1}{{2\sqrt 2 }}\exp \left( {\frac{{i\pi }}{4}} \right)} \right] - 18\,{\rm Li}_5 \left[ {\frac{1}{2}} \right] - \frac{9}{{16}}{\rm Li}_5 \left[ { - \frac{1}{4}} \right]\,.  
\end{split}
\end{equation}

Evaluating the identity at $(-1,i)$ and taking the real part gives

\begin{equation}\label{equ.deg5c}
\begin{split}
&\frac{{279}}{8}\,\zeta (5) - \frac{{977}}{{6144}}\,\pi ^4 \log 2 + \frac{{97}}{{768}}\,\pi ^2 \log ^3 2 - \frac{{15}}{{128}}\,\log ^5 2 \\ 
&\quad= 2\,{\mathop{\rm Re}\nolimits}\, {\rm Li}_5 \left[ {\frac{1}{{4\sqrt 2 }}\exp \left( {\frac{{i\pi }}{4}} \right)} \right] - 36\,{\mathop{\rm Re}\nolimits}\, {\rm Li}_5 \left[ {\frac{1}{{2\sqrt 2 }}\exp \left( {\frac{{i\pi }}{4}} \right)} \right] \\ 
&\qquad - 32\,{\mathop{\rm Re}\nolimits}\, {\rm Li}_5 \left[ {\frac{1}{{\sqrt 2 }}\exp \left( {\frac{{3i\pi }}{4}} \right)} \right] + 37\,{\rm Li}_5 \left[ {\frac{1}{2}} \right] + \frac{7}{8}\,{\rm Li}_5 \left[ { - \frac{1}{4}} \right]\,.
\end{split}
\end{equation}

Broadhurst proved (Eq.~68 of \cite{broadhurst98} written out) that

\begin{equation}\label{equ.deg5d}
\begin{split}
\frac{{31}}{{32}}\,\zeta (5) - \frac{{343}}{{99360}}\,\pi ^4 \log 2 + \frac{5}{{2484}}\,\pi ^2 \log ^3 2 - \frac{2}{{1035}}\,\log ^5 2 \\ 
  = \frac{{128}}{{69}}\,{\mathop{\rm Re}\nolimits}\, {\rm Li}_5 \left[ {\frac{1}{{\sqrt 2 }}\,\exp \left( {\frac{{i\pi }}{4}} \right)} \right] - \frac{{20}}{{69}}\,{\rm Li}_5 \left[ {\frac{1}{2}} \right] \,.
\end{split}
\end{equation}

\subsection{Base $2^{60}$ Binary BBP-type Formulas}

Solving Eqs.~\eqref{equ.deg5a}, \eqref{equ.deg5b}, \eqref{equ.deg5c} and \eqref{equ.deg5d} simultaneously for $\zeta (5)$, $\pi ^4 \log 2$, $\pi ^2 \log ^3 2$ and $\log ^5 2$ we find

\begin{equation}\label{equ.m3lvx5m}
\begin{split}
\zeta (5) &= \frac{{1317888}}{{1457}}\,{\mathop{\rm Re}\nolimits}\, {\rm Li}_5 \left[ {\frac{1}{{\sqrt 2 }}\exp \left( {\frac{{i\pi }}{4}} \right)} \right] - \frac{{377856}}{{62651}}\,{\mathop{\rm Re}\nolimits}\, {\rm Li}_5 \left[ {\frac{1}{{4\sqrt 2 }}\exp \left( {\frac{{i\pi }}{4}} \right)} \right] \\ 
  &\quad+ \frac{{56097792}}{{62651}}\,{\mathop{\rm Re}\nolimits}\, {\rm Li}_5 \left[ {\frac{1}{{\sqrt 2 }}\exp \left( {\frac{{3i\pi }}{4}} \right)} \right] + \frac{{1240064}}{{62651}}\,{\mathop{\rm Re}\nolimits}\, {\rm Li}_5 \left[ {\frac{1}{{2\sqrt 2 }}\exp \left( {\frac{{i\pi }}{4}} \right)} \right] \\ 
 &\quad - \frac{{929664}}{{62651}}{\rm Li}_5 \left[ {\frac{1}{2}} \right] + \frac{{644112}}{{62651}}{\rm Li}_5 \left[ {\frac{1}{4}} \right] \\ 
  &\quad- \frac{{63616}}{{62651}}{\rm Li}_5 \left[ { - \frac{1}{8}} \right] + \frac{{616752}}{{62651}}\,{\rm Li}_5 \left[ { - \frac{1}{4}} \right]\,,
\end{split}
\end{equation}

\begin{equation}\label{equ.fz69g9t}
\begin{split}
 \pi ^4 \log 2 &= \frac{{18593280}}{{47}}\,{\mathop{\rm Re}\nolimits}\, {\rm Li}_5 \left[ {\frac{1}{{\sqrt 2 }}\exp \left( {\frac{{i\pi }}{4}} \right)} \right] - \frac{{5294592}}{{2021}}\,{\mathop{\rm Re}\nolimits}\, {\rm Li}_5 \left[ {\frac{1}{{4\sqrt 2 }}\exp \left( {\frac{{i\pi }}{4}} \right)} \right] \\ 
 &\quad + \frac{{794230272}}{{2021}}\,{\mathop{\rm Re}\nolimits}\, {\rm Li}_5 \left[ {\frac{1}{{\sqrt 2 }}\exp \left( {\frac{{3i\pi }}{4}} \right)} \right] + \frac{{16467456}}{{2021}}\,{\mathop{\rm Re}\nolimits}\, {\rm Li}_5 \left[ {\frac{1}{{2\sqrt 2 }}\exp \left( {\frac{{i\pi }}{4}} \right)} \right] \\ 
 &\quad - \frac{{11408256}}{{2021}}{\rm Li}_5 \left[ {\frac{1}{2}} \right] + \frac{{9121248}}{{2021}}{\rm Li}_5 \left[ {\frac{1}{4}} \right] \\ 
 & \quad- \frac{{900864}}{{2021}}{\rm Li}_5 \left[ { - \frac{1}{8}} \right] + \frac{{8769816}}{{2021}}\,{\rm Li}_5 \left[ { - \frac{1}{4}} \right]\,,
\end{split}
\end{equation}

\begin{equation}\label{equ.nlux2ag}
\begin{split}
 \pi ^2 \log ^3 2 &= \frac{{17440704}}{{47}}\,{\mathop{\rm Re}\nolimits}\, {\rm Li}_5 \left[ {\frac{1}{{\sqrt 2 }}\exp \left( {\frac{{i\pi }}{4}} \right)} \right] - \frac{{4835904}}{{2021}}\,{\mathop{\rm Re}\nolimits}\, {\rm Li}_5 \left[ {\frac{1}{{4\sqrt 2 }}\exp \left( {\frac{{i\pi }}{4}} \right)} \right] \\ 
 &\quad + \frac{{747218112}}{{2021}}\,{\mathop{\rm Re}\nolimits}\, {\rm Li}_5 \left[ {\frac{1}{{\sqrt 2 }}\exp \left( {\frac{{3i\pi }}{4}} \right)} \right] + \frac{{12619200}}{{2021}}\,{\mathop{\rm Re}\nolimits}\, {\rm Li}_5 \left[ {\frac{1}{{2\sqrt 2 }}\exp \left( {\frac{{i\pi }}{4}} \right)} \right] \\ 
  &\quad- \frac{{7432176}}{{2021}}{\rm Li}_5 \left[ {\frac{1}{2}} \right] + \frac{{8731962}}{{2021}}{\rm Li}_5 \left[ {\frac{1}{4}} \right] \\ 
  &\quad- \frac{{862416}}{{2021}}{\rm Li}_5 \left[ { - \frac{1}{8}} \right] + \frac{{8350599}}{{2021}}\,{\rm Li}_5 \left[ { - \frac{1}{4}} \right]
\end{split}
\end{equation}

and

\begin{equation}\label{equ.mspv1c6}
\begin{split}
 \log ^5 2 &= \frac{{6218880}}{{47}}\,{\mathop{\rm Re}\nolimits}\, {\rm Li}_5 \left[ {\frac{1}{{\sqrt 2 }}\exp \left( {\frac{{i\pi }}{4}} \right)} \right] - \frac{{1689472}}{{2021}}\,{\mathop{\rm Re}\nolimits}\, {\rm Li}_5 \left[ {\frac{1}{{4\sqrt 2 }}\exp \left( {\frac{{i\pi }}{4}} \right)} \right] \\ 
  &\quad+ \frac{{266699392}}{{2021}}\,{\mathop{\rm Re}\nolimits}\, {\rm Li}_5 \left[ {\frac{1}{{\sqrt 2 }}\exp \left( {\frac{{3i\pi }}{4}} \right)} \right] + \frac{{3780736}}{{2021}}\,{\mathop{\rm Re}\nolimits}\, {\rm Li}_5 \left[ {\frac{1}{{2\sqrt 2 }}\exp \left( {\frac{{i\pi }}{4}} \right)} \right] \\ 
 &\quad - \frac{{2092736}}{{2021}}{\rm Li}_5 \left[ {\frac{1}{2}} \right] + \frac{{3217563}}{{2021}}{\rm Li}_5 \left[ {\frac{1}{4}} \right] \\ 
  &\quad- \frac{{317784}}{{2021}}{\rm Li}_5 \left[ { - \frac{1}{8}} \right] + \frac{{3005666}}{{2021}}\,{\rm Li}_5 \left[ { - \frac{1}{4}} \right]\,.
\end{split}
\end{equation}

Identities \eqref{equ.m3lvx5m}, \eqref{equ.fz69g9t}, \eqref{equ.nlux2ag} and \eqref{equ.mspv1c6} facilitate the derivation of base $2^{60}$, length $120$ BBP-type formulas for the respective polylogarithm constants. The explicit formulas are as listed in the BBP Compendium.

\section{Zero Relations}

BBP zero relations are BBP-type formulas that evaluate to zero. Considering that BBP-type formulas are usually discovered through computer searches, the need to study zero relations is aptly set forth in the BBP Compendium~\cite{bailey09}:

\begin{quote}
Knowledge of these zero relations is essential for finding formulas using integer relation programs (such as PSLQ). This is because unless these zero relations are excluded from the search for a conjectured BBP-type formula, the search may only recover a zero relation.
\end{quote}

\subsection{Degree $2$ Zero Relations}

\subsubsection{Base $2^{12}$ Relations}

Eliminating $\pi^2$ between Eqs.~\eqref{equ.whrutog} and \eqref{equ.dmg73yi}, we find the identity

\begin{eqnarray}
 0 &=& 72\,{\rm Re\,Li}_2 \left[ {\frac{1}{2}\exp \left( {\frac{{i\pi }}{2}} \right)} \right] - 72\,{\rm Re\,Li}_2 \left[ {\frac{1}{{\sqrt 2 }}\exp \left( {\frac{{3i\pi }}{4}} \right)} \right] \nonumber\\ 
&&  - 24\,{\rm Re\,Li}_2 \left[ {\frac{1}{{2\sqrt 2 }}\exp \left( {\frac{{i\pi }}{4}} \right)} \right] - 180\,{\rm Re\,Li}_2 \left[ {\frac{1}{2}\exp \left( {\frac{{i\pi }}{3}} \right)} \right]\nonumber \\ 
&&  + 45\,{\rm Li}_2 \left[ {\frac{1}{4}} \right]\,,
\end{eqnarray}

which leads immediately to the zero relation

\begin{eqnarray}\label{equ.hbrwred}
0&=&P(2,2^{12},24,(2^{11},-5\cdot2^{11},-2^{12},3\cdot2^{12},-2^{9},5\cdot2^{10},2^{8},\nonumber\\
&& 3^2\cdot2^{9},2^{9},-5\cdot2^{7},-2^{6},0,-2^{5},-5\cdot2^{5},2^{6},3^2\cdot2^{5},2^{3},\nonumber\\
&& 5\cdot2^{4},-2^{2},3\cdot2^{4},-2^{3},-5\cdot2,1,0))\,.
\end{eqnarray}

Eliminating $\log^2 2$ between Eq.~\eqref{equ.nptw12y} and \eqref{equ.c5ut706} gives the identity

\begin{eqnarray}
 0 &=& 2\,{\rm Re\,Li}_2 \left[ {\frac{1}{{\sqrt 2 }}\exp \left( {\frac{{3i\pi }}{4}} \right)} \right] + 2\,{\rm ReLi}_2 \left[ {\frac{1}{{2\sqrt 2 }}\exp \left( {\frac{{i\pi }}{4}} \right)} \right] \nonumber\\ 
&&  + {\rm Li}_2 \left[ { - \frac{1}{8}} \right] - 2\,{\rm Li}_2 \left[ { - \frac{1}{2}} \right] - 2\,{\rm Li}_2 \left[ {\frac{1}{4}} \right] - {\rm Li}_2 \left[ { - \frac{1}{4}} \right]\,, 
\end{eqnarray}

which produces the zero relation

\begin{eqnarray}\label{equ.m1x0ete}
0&=&P(2,2^{12},24,(2^{11},-2^{13},-5\cdot2^{11},13\cdot2^{10},-2^{9},7\cdot2^{10},2^{8},3^3\cdot2^{8},\nonumber\\
&& 5\cdot2^{8},-2^{9},-2^{6},2^{8},-2^{5},-2^{7},5\cdot2^{5},3^3\cdot2^{4},2^{3},7\cdot2^{4},\nonumber\\
&& -2^{2},13\cdot2^{2},-5\cdot2^{2},-2^{3},1,0))\,.
\end{eqnarray}

\subsubsection{A Base $2^{60}$ Relation}

Eliminating $\pi\log 2$ between Eq.~\eqref{equ.jkxe2jb} and Eq.~\eqref{equ.ibwzeu0} gives the identity

\begin{eqnarray}\label{equ.w33gn39}
0 &=& 5\,{\rm Im\,Li}_2 \left[ {\frac{1}{2}\exp \left( {\frac{{i\pi }}{2}} \right)} \right] - 4\,{\rm Im\,Li}_2 \left[ {\frac{1}{{\sqrt 2 }}\exp \left( {\frac{{i3\pi }}{4}} \right)} \right] \nonumber\\ 
&& - 2\,{\rm Im\,Li}_2 \left[ {\frac{1}{{2\sqrt 2 }}\exp \left( {\frac{{i\pi }}{4}} \right)} \right] - 2\,{\rm Im\,Li}_2 \left[ {\frac{1}{{4\sqrt 2 }}\exp \left( {\frac{{i\pi }}{4}} \right)} \right]\,, 
\end{eqnarray}

from which we get, immediately, the binary BBP-type zero relation

\begin{eqnarray}
&& 0=\nonumber\\
&& P(2,2^{60},120,(-2^{59},3 \cdot 2^{60},-11 \cdot 2^{57},0,-23 \cdot 2^{56},-3 \cdot 7 \cdot 2^{56},2^{56},0,-11 \cdot 2^{54},-13 \cdot 2^{54},\nonumber\\
&& -2^{54},0,2^{53},-3 \cdot 2^{54},-7 \cdot 2^{52},0,-2^{51},3 \cdot 7 \cdot 2^{50},-2^{50},0,11 \cdot 2^{48},-3 \cdot 2^{50},2^{48},0,\nonumber\\
&& 23 \cdot 2^{46},3 \cdot 2^{48},-11 \cdot 2^{45},0,2^{45},2^{46},2^{44},0,-11 \cdot 2^{42},3 \cdot 2^{44},23 \cdot 2^{41},0,2^{41},-3 \cdot 2^{42},\nonumber\\
&& 11 \cdot 2^{39},0,-2^{39},3 \cdot 7 \cdot 2^{38},-2^{38},0,-7 \cdot 2^{37},-3 \cdot 2^{38},2^{36},0,-2^{35},-13 \cdot 2^{34},-11 \cdot 2^{33},0,\nonumber\\
&& 2^{33},-3 \cdot 7 \cdot 2^{32},-23 \cdot 2^{31},0,-11 \cdot 2^{30},3 \cdot 2^{32},-2^{30},0,2^{29},-3 \cdot 2^{30},11 \cdot 2^{27},0,23 \cdot 2^{26},\nonumber\\
&& 3 \cdot 7 \cdot 2^{26},-2^{26},0,11 \cdot 2^{24},13 \cdot 2^{24},2^{24},0,-2^{23},3 \cdot 2^{24},7 \cdot 2^{22},0,2^{21},-3 \cdot 7 \cdot 2^{20},2^{20},\nonumber\\
&& 0,-11 \cdot 2^{18},3 \cdot 2^{20},-2^{18},0,-23 \cdot 2^{16},-3 \cdot 2^{18},11 \cdot 2^{15},0,-2^{15},-2^{16},-2^{14},0,11 \cdot 2^{12},\nonumber\\
&& -3 \cdot 2^{14},-23 \cdot 2^{11},0,-2^{11},3 \cdot 2^{12},-11 \cdot 2^{9},0,2^{9},-3 \cdot 7 \cdot 2^{8},2^{8},0,7 \cdot 2^{7},3 \cdot 2^{8},-2^{6},0,\nonumber\\
&& 2^{5},13 \cdot 2^{4},11 \cdot 2^{3},0,-2^{3},3 \cdot 7 \cdot 2^{2},23 \cdot 2,0,11,-3 \cdot 2^{2},1,0))\,.
\end{eqnarray}

\subsection{Degree $3$ Zero Relations}

\subsubsection{Base $2^{12}$ Zero Relation}

Eliminating $\pi ^2 \log 2$ and $\log^3 2$ between Eqs.~\eqref{equ.reeuf8k}, \eqref{equ.hsmmff3} and \eqref{equ.mqgomtf}, and subtracting the resulting identity from Eq.~\eqref{equ.ls90hkt} gives the BBP-ready identity

\begin{eqnarray}
 0 &=&  - \frac{{88}}{{21}}\,{\rm Li}_3 \left[ {\frac{1}{2}} \right] - \frac{{132}}{{133}}\,{\rm Li}_3 \left[ {\frac{1}{4}} \right] + \frac{{88}}{{399}}\,{\rm Li}_3 \left[ { - \frac{1}{8}} \right] - \frac{{1584}}{{133}}\,{\rm Re\,Li}_3 \left[ {\frac{1}{2}\exp \left( {\frac{{i\pi }}{2}} \right)} \right] \nonumber\\ 
&&  + \frac{{704}}{{133}}\,{\rm Re\,Li}_3 \left[ {\frac{1}{{2\sqrt 2 }}\exp \left( {\frac{{i\pi }}{4}} \right)} \right] + \frac{{704}}{{399}}\,{\rm Re\,Li}_3 \left[ {\frac{1}{{\sqrt 2 }}\exp \left( {\frac{{i\pi }}{4}} \right)} \right]\,, 
 \end{eqnarray}

which yields the only {\em expected} degree $3$ base $2^{12}$ length $24$ zero relation listed in the BBP Compendium and which was originally discovered by Bailey~\cite{bailey09}, using his PSLQ program.

\subsubsection{Base $2^{60}$ Zero Relations}

Solving Eqs.~\eqref{equ.iex39kd}, \eqref{equ.w2k1gjt} and \eqref{equ.mqgomtf} for $\zeta(3)$ and subtracting Eq.~\eqref{equ.fgc5kj5} gives the identity

\begin{eqnarray}\label{equ.pblqhqk}
 0 &=& \frac{{50}}{{21}}\,{\rm Li}_3 \left[ {\frac{1}{2}} \right] - \frac{5}{7}\,{\rm Li}_3 \left[ {\frac{1}{4}} \right] + \frac{{10}}{{63}}\,{\rm Li}_3 \left[ { - \frac{1}{8}} \right] - \frac{{319}}{{63}}\,{\rm Re\,Li}_3 \left[ {\frac{1}{2}\exp \left( {\frac{{i\pi }}{2}} \right)} \right]\nonumber \\ 
&&  - \frac{{236}}{{63}}\,{\rm Re\,Li}_3 \left[ {\frac{1}{{2\sqrt 2 }}\exp \left( {\frac{{i\pi }}{4}} \right)} \right] - \frac{8}{9}\,{\rm Re\,Li}_3 \left[ {\frac{1}{{\sqrt 2 }}\exp \left( {\frac{{i\pi }}{4}} \right)} \right]\nonumber \\ 
&&  + \frac{{68}}{{63}}\,{\rm Re\,Li}_3 \left[ {\frac{1}{{4\sqrt 2 }}\exp \left( {\frac{{i\pi }}{4}} \right)} \right]\,.
\end{eqnarray}

Solving Eqs.~\eqref{equ.hsmmff3}, \eqref{equ.w2k1gjt} and \eqref{equ.mqgomtf} for $\zeta(3)$ and subtracting Eq.~\eqref{equ.fgc5kj5} gives the identity

\begin{eqnarray}\label{equ.v0bcwbn}
 0 &=&  - \frac{{30}}{{119}}\,{\rm Li}_3 \left[ {\frac{1}{4}} \right] + \frac{{20}}{{357}}\,{\rm Li}_3 \left[ { - \frac{1}{8}} \right] - \frac{7}{3}\,{\rm Re\, Li}_3 \left[ {\frac{1}{2}\exp \left( {\frac{{i\pi }}{2}} \right)} \right] \nonumber\\ 
&&  - \frac{{52}}{{357}}\,{\rm Re\, Li}_3 \left[ {\frac{1}{{2\sqrt 2 }}\exp \left( {\frac{{i\pi }}{4}} \right)} \right] + \frac{8}{{357}}\,{\rm Re\, Li}_3 \left[ {\frac{1}{{\sqrt 2 }}\exp \left( {\frac{{i\pi }}{4}} \right)} \right]\nonumber \\ 
&&  + \frac{{76}}{{357}}\,{\rm Re\, Li}_3 \left[ {\frac{1}{{4\sqrt 2 }}\exp \left( {\frac{{i\pi }}{4}} \right)} \right]\,.
\end{eqnarray}

The explicit base $2^{60}$ zero relations resulting from the identities \eqref{equ.pblqhqk} and \eqref{equ.v0bcwbn} are as listed in the BBP Compendium.

\subsection{Degree $4$ Zero Relation}

\subsubsection{A Base $2^{60}$ Relation}

Subtracting Eq.~\eqref{equ.b2nxb0w} from Eq.~\eqref{equ.b3ywont}, we obtain the identity

\begin{equation}\label{equ.p6i7d5p}
\begin{split}
0 &=  - \frac{{398592}}{{2911}}{\mathop{\rm Re}\nolimits} {\rm Li}_4 \left[ {\frac{1}{{2\sqrt 2 }}\exp \left( {\frac{{i\pi }}{4}} \right)} \right] + \frac{{7627392}}{{20377}}{\mathop{\rm Re}\nolimits} {\rm Li}_4 \left[ {\frac{1}{{\sqrt 2 }}\exp \left( {\frac{{i\pi }}{4}} \right)} \right] \\ 
 &\qquad + \frac{{1150848}}{{2911}}{\mathop{\rm Re}\nolimits} {\rm Li}_4 \left[ {\frac{1}{{\sqrt 2 }}\exp \left( {\frac{{3i\pi }}{4}} \right)} \right] + \frac{{34560}}{{497}}{\mathop{\rm Re}\nolimits} {\rm Li}_4 \left[ {\frac{1}{{4\sqrt 2 }}\exp \left( {\frac{{i\pi }}{4}} \right)} \right] \\ 
  &\qquad- \frac{{38797488}}{{20377}}{\mathop{\rm Re}\nolimits} {\rm Li}_4 \left[ {\frac{1}{2}\exp \left( {\frac{{i\pi }}{2}} \right)} \right] + \frac{{1394496}}{{20377}}{\rm Li}_4 \left[ {\frac{1}{2}} \right] \\ 
  &\qquad- \frac{{3888}}{{41}}{\rm Li}_4 \left[ {\frac{1}{4}} \right] + \frac{{576}}{{41}}{\rm Li}_4 \left[ { - \frac{1}{8}} \right] + \frac{{26136}}{{2911}}{\rm Li}_4 \left[ { - \frac{1}{4}} \right] \,.
\end{split}
\end{equation}

As usual, writing each component of Eq.~\eqref{equ.p6i7d5p} as a base~$2^{60}$, length~$120$ binary BBP-type formula and forming the indicated combination, we obtain the only degree $4$ binary zero relation, given explicitly in the BBP Compendium.

\section{Conclusion}

Using a clear and straightforward approach, we have obtained and proved interesting new binary digit extraction BBP-type formulas for polylogarithm constants. Some known results were also rediscovered in a clearer and more elegant manner. Experimentally discovered binary BBP-type formulas are also proved.

\section{Acknowledgement}
The author enjoyed interesting correspondence with Dr.\ David H.\ Bailey and Jaume O.\ Lafont. He is also grateful to the reviewer for a detailed excellent review with very useful comments and suggestions which led to a significantly improved manuscript. The author received a lot of support, encouragement and goodwill from Prof.\ H.\ W.\ Lenstra.

\end{document}